\documentclass[a4paper]{amsart}
\usepackage{orcidlink}
 \usepackage[foot]{amsaddr}

\usepackage[backend=biber,style=numeric, sorting=nyt, maxnames=10]{biblatex}
\usepackage{graphicx}
\graphicspath{{Figures/}}

\usepackage{amsmath}
\usepackage{amssymb}
\usepackage{float}
\usepackage{csquotes}

\usepackage{todonotes}

\usepackage{algorithm}
\usepackage{algpseudocode}

\usepackage{subcaption}

\usepackage{tikz}

\usetikzlibrary{arrows.meta, decorations.pathreplacing}

\usepackage{booktabs}

\newcommand{\R}{\mathbb{R}}

\title[Adjoint-based Perfusion Estimation from DCE-US]{Adjoint-based Perfusion Estimation from Dynamic Contrast-Enhanced Ultrasound: Advection-Diffusion and Two-Compartment Models}

\author[Externbrink et al.]{Sophie Externbrink$^1$\orcidlink{0000-0000-0000-0000}}
\author[]{Ahmed El Kaffas$^2$\orcidlink{0000-0001-5296-0247}}
\author[]{Dimitre Hristov$^3$\orcidlink{0000-0002-2462-7727}}
\author[]{Sebastian Götschel$^1$\orcidlink{0000-0003-0287-2120}}

\address{$^1$Chair Computational Mathematics, Institute of Mathematics, Hamburg University of Technology, 21073 Hamburg, Germany}

\address{$^2$Department of Radiology, University of California, San Diego, CA, USA}

\address{$^3$Radiation Oncology, Stanford University Medical Center, Stanford, CA, USA}

\email{sophie.externbrink@tuhh.de}
\email{aelkaffas@health.ucsd.edu}
\email{dhristov@stanford.edu}
\email{sebastian.goetschel@tuhh.de}

\begin{document}

\begin{abstract}
	Tumor perfusion and vascular properties are important determinants of a cancer's response to therapy. 
	In this paper, we discuss the estimation of spatially varying blood flow velocities and perfusion parameters from time-resolved contrast agent concentration data. We compare a standard parabolic advection-diffusion model against a  two-compartment model governed by a coupled system of hyperbolic advection-reaction equations, which is physiologically more sound. To address the inherent ill-posedness of this parameter identification problem, we employ Tikhonov regularization and derive continuous adjoint equations necessary for efficient, gradient-based minimization. We discuss the numerical discretization of the state and adjoint systems using state-of-the-art schemes, and demonstrate the efficacy of the proposed reconstruction algorithms through numerical experiments on synthetic data and \emph{in vivo} dynamic contrast-enhanced ultrasound measurements. 
\end{abstract}

\maketitle

\section{Introduction}
Liver cancer is a major cause of cancer-related mortality due to its poor prognosis, despite its relatively low incidence. 
In the European Union, approximately 62,000 new cases were diagnosed in 2022, and around 54,000 deaths from liver cancer are estimated for 2022~\cite{ECIS}.
Liver cancer occurrence is heavily linked with obesity, type-2 diabetes, lack of movement, alcohol consumption
and hepatic viral infections. Current studies suggest that these factors will further increase liver cancer
cases world-wide by more than 55\% by 2040 and a doubling of nonalcoholic fatty liver caused cancers in France, USA and China by 2030~\cite{RUMGAY20221598,HuangDanielQ.2021GeoN}. 
Without appropriate treatment tailored to the tumor type, the average survival rate is only 8.7 months \cite{LiverCancerKills}. 
The tumor vascular system, often described as the ``fingerprints" of cancer, exhibits significant variability across cancer types
and patients, profoundly influencing the response to various therapies, including classical radiation and chemotherapy as well 
as newer immunotherapies \cite{FingerprintCancer, TumourVascularImportant}. Therefore, optimizing treatment on an individual 
patient basis is crucial. It makes early adaptations or change in therapy possible, when treatment proves ineffective
and thus, reduces the often severe side effects of cancer treatment \cite{PatientCare}. However, evaluating treatment 
efficiency is challenging. Imaging modalities such as computed tomography (CT), magnetic resonance imaging (MRI), 
and positron emission tomography (PET) allow direct visualization of tumor regions to assess response but are costly, 
involve radiation exposure (CT/PET), and are not suitable for frequent bedside monitoring \cite{Scans}. 
These limitations restrict their use during ongoing treatment. Alternatively, tumor perfusion 
and vascular characteristics can serve as indirect indicators of therapeutic effect, since tumor growth and survival depend 
heavily on blood supply. 

Dynamic contrast-enhanced ultrasound (DCE-US) provides a cost-effective, radiation-free approach to 
visualize tumor blood flow and perfusion \cite{bloodflow}. DCE-US is a functional imaging technique that uses contrast agents (microbubbles) injected into the blood, and records the movement of such a contrast bolus through an organ or tumor. Besides heuristic techniques to quantify perfusion from such measurements, like considering the mean transit time, the movement of the contrast agent can be 
mathematically modeled using tracer kinetic models \cite{Tracer}. Employing such models enables reconstruction of tumor perfusion and vascular parameters on a patient-specific basis, facilitating improved individualised therapy.
Reconstructing these parameters involves solving an inverse problem: given tracer concentration data over time and a mathematical model, parameter functions are identified through fitting simulation and measurement data. 
Although often ill-posed, parameter identification problems have been extensively studied and successfully applied across 
diverse fields—including prostate cancer diagnosis \cite{prostateCancer}, ocean engineering \cite{oceanEngineering}, and 
environmental monitoring \cite{ExmapleAdvecDiff}. In particular, advection-diffusion models have been widely used to simulate 
transport phenomena, with early computational reconstructions of advective flows dating back to the early 2000s \cite{AdvecExmaple3}. 
Applications range from tracer dispersion \cite{AdvecExample} to the formation and evolution of orogenic topography \cite{AdvecExample4}.

Selecting an appropriate forward model is critical for obtaining physiologically meaningful results applicable to 
cancer therapy. While advection-diffusion tracer models are common, they often fail to fully capture the physical 
realities of tracer transport within organs, relying on region-of-interest frameworks that treat the system as 
isolated with a single arterial input function \cite{DiffusionExample1, AdvecDiffIntro}. Hristov et al.~\cite{bloodflow} applied an advection-diffusion model~\cite{Sourbron} to reconstruct parameters sensitive 
to vascular network changes correlated with histology, yet the diffusion parameter lacked clear physiological interpretation 
and showed only weak correlations. To overcome these limitations, 
Sourbron proposed alternative models including a two-compartment framework \cite{Sourbron}. 
This model separates tracer concentration into arterial and venous compartments, enabling the description of 
tracer transport via conversion from arterial to venous blood rather than pure diffusion. A space-dependent transfer 
function characterises this conversion for a more suitable modeling of bolus dispersion, and thereby improving physiological relevance.

In this paper we employ Sourbron’s two-compartment model as the basis for parameter identification, and compare its results to the advection-diffusion model. We reconstruct blood flow velocities and the conversion parameter as spatially varying functions from time-dependent concentration measurements. To solve the arising  partial differential equation (PDE)-constrained optimization problems, we derive adjoint equations for computationally efficient gradient computation. We provide a prototype implementation using state-of-the-art discretization techniques to solve the PDEs numerically. The reconstruction is tested on synthetic data as well as DCE-US measurement data. Both mathematical models are introduced in Section~\ref{sec:model}, the parameter reconstruction methods are discussed in Section~\ref{sec:reconstruction}, and numerical results presented in Section~\ref{sec:results}.

\section{Mathematical models: Tracer-kinetic field theory} \label{sec:model}

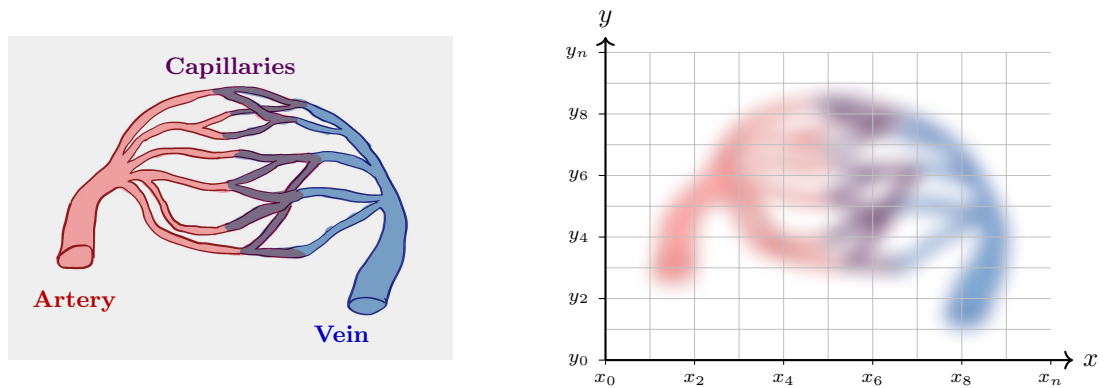
\begin{figure}[!htp]
  \centering
  \begin{subfigure}[c]{0.48\textwidth}
        \centering
        \begin{tikzpicture}

          \node[anchor=south west, inner sep=0] (img) 
            {\includegraphics[width=0.8\textwidth]{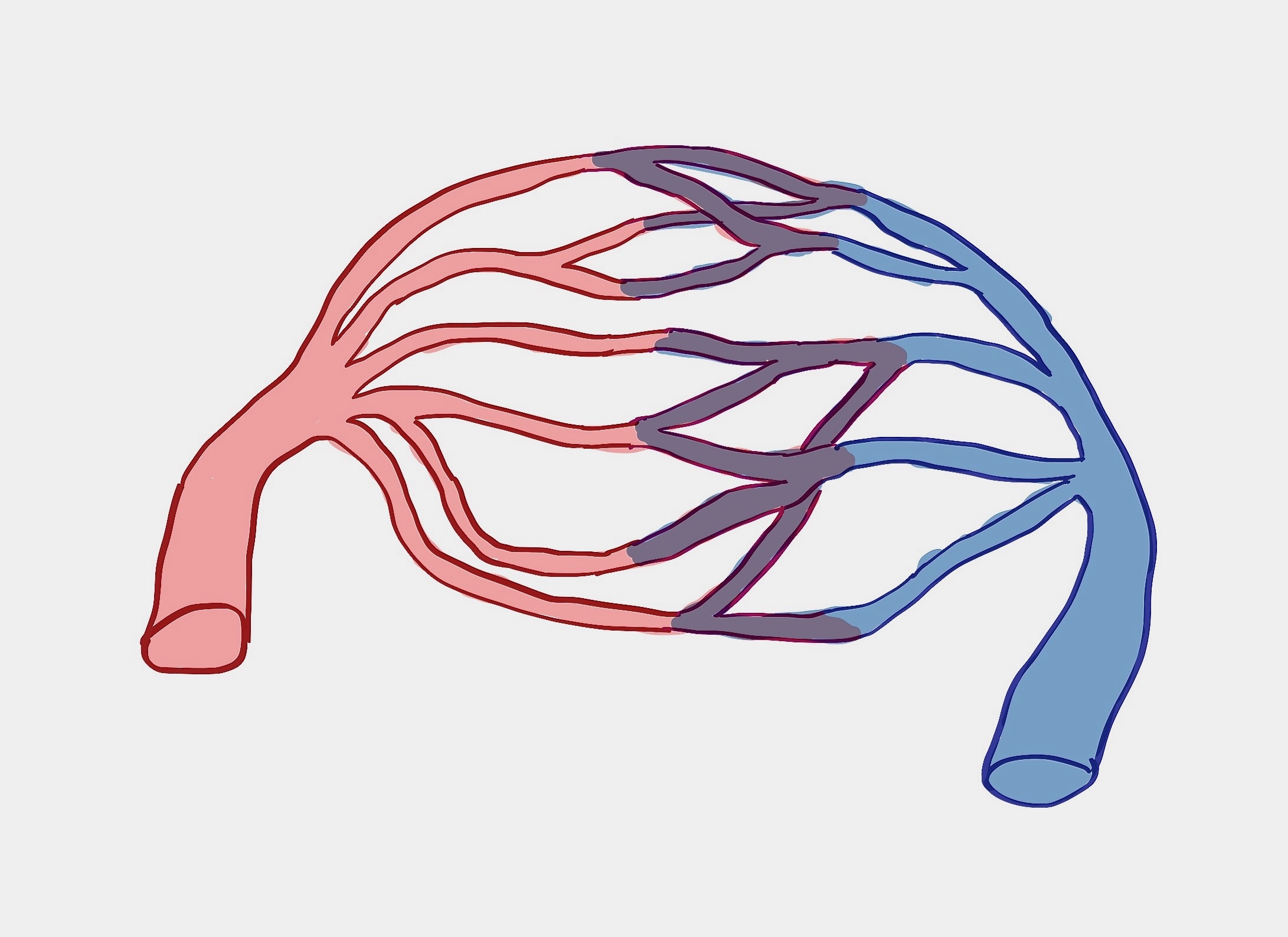}};

          \begin{scope}[x={(img.south east)}, y={(img.north west)}]

          \node[font=\bfseries\small, text=red!70!black]
            at (0.15,0.18) {Artery};

          \node[font=\bfseries\small, text=blue!70!black]
            at (0.75,0.08) {Vein};

          \node[font=\bfseries\small, text=violet!70!black]
            at (0.50,0.90) {Capillaries};

        \end{scope}
        \end{tikzpicture}
    \end{subfigure}
    \hfill
    \begin{subfigure}[c]{0.48\textwidth}
        \centering
        \begin{tikzpicture}
          \node[anchor=south west, inner sep=0] (img)
            {\includegraphics[width=0.8\textwidth]{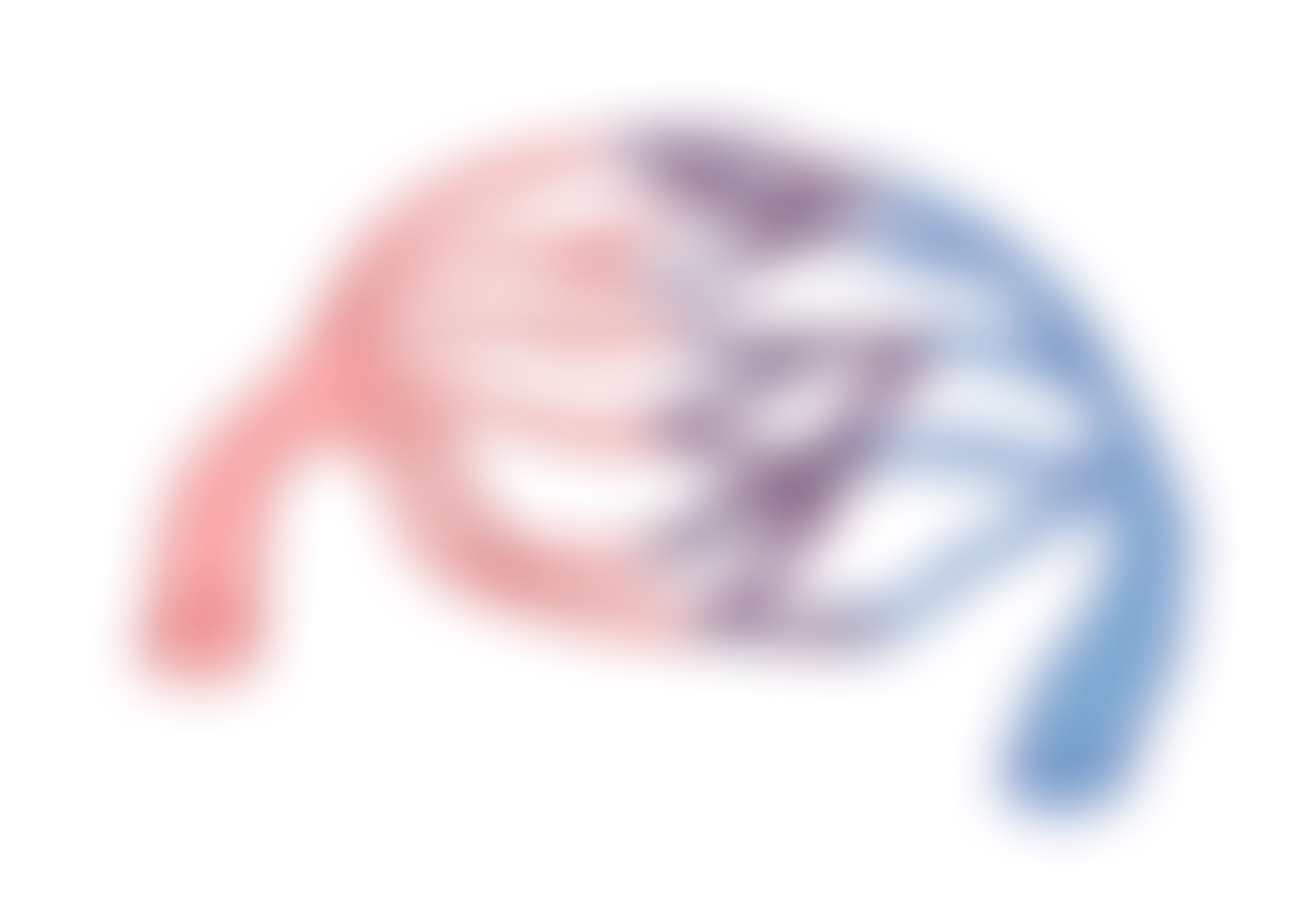}};

          \begin{scope}[x={(img.south east)}, y={(img.north west)}]

            \draw[
              step=0.1,
              gray!50,
              very thin
            ] (0,0) grid (1,1);

            \draw[->, thick] (0,0) -- (1.05,0) node[right] {$x$};
            \draw[->, thick] (0,0) -- (0,1.05) node[above] {$y$};
            
            \foreach \i/\k in {0/0,0.2/2,0.4/4,0.6/6, 0.8/8,1/n} {
            \draw (\i,0) -- (\i,-0.015);
            \node[below] at (\i,-0.015)
              {\scriptsize $x_{\k}$};
            }

            \foreach \i/\k in {0/0,0.2/2,0.4/4,0.6/6, 0.8/8,1/n} {
              \draw (0,\i) -- (-0.015,\i);
              \node[left] at (-0.015,\i)
                {\scriptsize $y_{\k}$};
  }

          \end{scope}
          \end{tikzpicture}
    \end{subfigure}
    \caption{Left: Schematic representation of a vascular system, vessel color indicating presence of oxygenated (red) and de-oxygenated (blue) blood. Right: The finite resolution of the imaging system (ultrasound) eliminates small vessel details.}
    \label{fig:liverSimple}
\end{figure}

In this section we briefly summarize two mathematical models to describe the propagation of tracer concentrations in blood flow. Contrast agent is washed into the organ or tumor through an arterial inlet, spreads out through the capillaries, where blood gets de-oxygenated, and then is washed out through a venous outlet, as illustrated in Figure \ref{fig:liverSimple}. Note that the ultrasound imaging system typically will not be able to resolve small vessel details.

\subsection{Advection-diffusion model}

\begin{figure}[!htp]
  \centering
  \includegraphics[width = 0.49\textwidth]{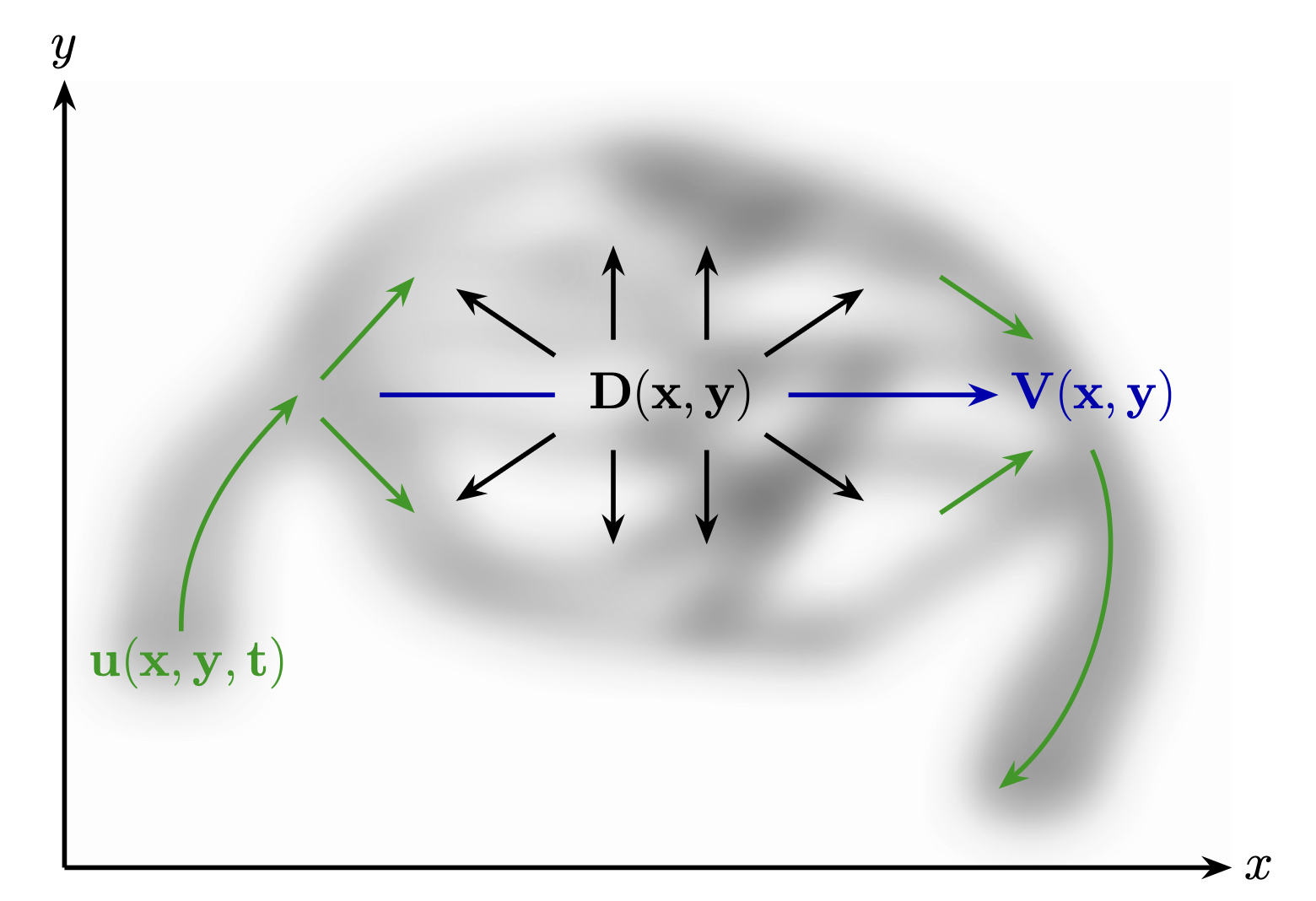}  \includegraphics[width=0.49\textwidth]{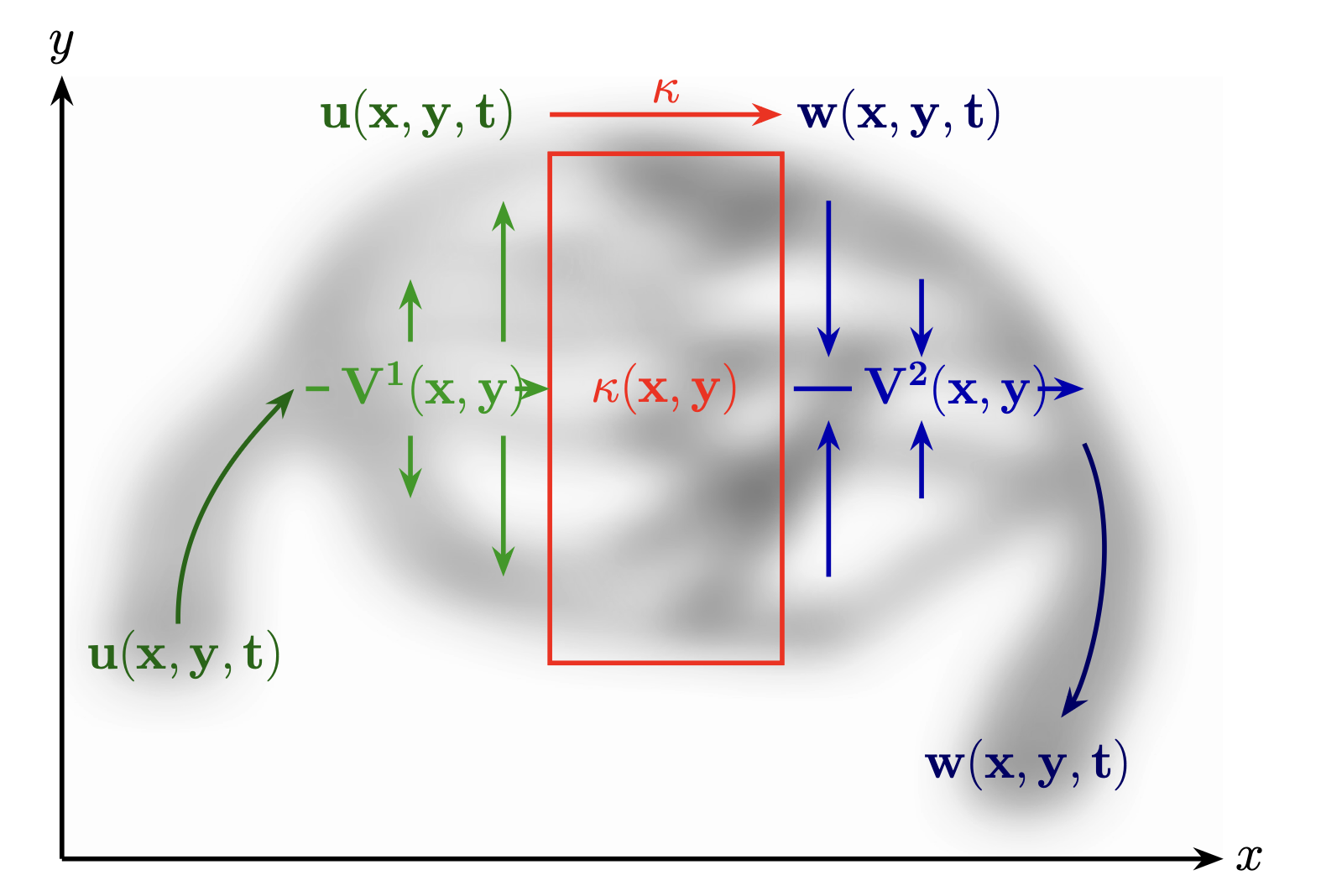}
  \caption{Sketch of transportation of tracer in blood flow through an organ using advection-diffusion (left) and the two-compartment model (right).}
  \label{fig:sketch_diff}
\end{figure}

When modeling tracer in blood flow of an organ we need to be able to describe an 
in- and outflow via the main artery, respectively vein. Furthermore, upon entering the organ the tracer spreads out 
through a lot of smaller vessels and capillaries. As in measurements and simulations, these small blood vessels often cannot be resolved due to technical and computational restrictions, the behavior often is approximated by diffusion. 
The general behavior can be seen in Figure~\ref{fig:sketch_diff}~(left). For the tracer concentration $u(x,t)$, the advection-diffusion model on a space-time domain $Q = \Omega \times (0,T)$ is given by
\begin{equation}\label{eq:advdiff}
u_{t} + \nabla \cdot (V \, u) - \nabla \cdot (D \, \nabla u) = 0 \text{ on } Q,
\end{equation}
together with suitable initial and boundary conditions. Here, $\Omega \subset \R^d$ denotes the spatial domain, typically a rectangular domain given by the imaging modality, and composed of pixels ($d=2$) or voxels ($d=3$).
$V : \Omega \to \R^d$ is the velocity field, i.e., the heart-beat averaged blood flow velocity, and 
$D$ the diffusion $D: \Omega \rightarrow \mathbb{R}$, which we assume to be isotropic~\cite{Sourbron}. Inflow from and outflow to outside the imaging domain can be modeled by including suitable boundary conditions or source term. However, as  these additional functions would need to be reconstructed from measurements together with velocities and diffusivity, this increases the computational complexity significantly, and potentially would require prior knowledge from other imaging modalities. For the remainder of this paper, we thus assume that the contrast agent bolus is fully contained inside the computational domain,~i.e., consider time intervals not including the inflow and outflow of the bolus.

\subsection{Two-compartment model}

While initial movement through capillaries can be modeled via a diffusive part, modeling the tracer exiting the capillaries,~i.e., merging into the vein, requires modifications to avoid negative diffusivities. Also while moderate correlation is shown in~\cite{bloodflow}, there is no direct link to the medically relevant parameter of perfusion. 
As a remedy, a two-compartment model has been proposed~\cite{Sourbron}.
%
%
%
The model consists of two coupled advection equations, where the two
 components, $u$ and $w$,  represent tracer concentrations in the oxygenated (arterial) and de-oxygenated (venous) blood, respectively.
To get a realistic model of an organ, we make sure that there is only oxygenated blood in the inflow region and only de-oxygenated blood in 
the outflow region. The oxygenated blood is converted into de-oxygenated blood using a space-dependent rate (transfer coefficient function) $\kappa : \Omega \to \R$. 
In addition, both components can move with different velocities $V^1$ and $V^2$ for the arterial and venous blood flow, 
respectively, see Figure~\ref{fig:sketch_diff} (right). When modeling the behavior of blood flow through an organ the velocity of the
oxygenated blood should be fast at the inflow boundary and reduced to zero at the outflow boundary. Opposite to this, the
de-oxygenated blood should be slow at the inflow boundary and fast at the outflow boundary, see Table \ref{fig:bloodFlowTable}.
%
\begin{table}
    \centering
    \begin{tabular}{ p{3cm}  p{3cm}} 
        \toprule
        Vessel Type & Mean Flow Velocity \\ 
        \midrule
        Aorta      & $30 \, cm/s$ \\
        Arteries   & $30 \, cm/s$ - $20 \, mm/s$ \\
        Arterioles & $20 \, mm/s$ \\
        \midrule
        Capillaries & $0.3 \, mm/s$ \\
        \midrule
        Venules   & $3 \, mm/s$ \\
        Veins     & $3 \, mm/s$ - $10 \, mm/s$ \\
        Vena Cava & $10 \, mm/s$ \\
        \bottomrule
    \end{tabular}
    \caption{Average haemodynamic parameters of the human circulation \cite{BloosflowBook}. Only in the capillaries, the transfer coefficient $\kappa$ of the two-compartment model in nonzero.}  
    \label{fig:bloodFlowTable}
\end{table}

Altogether, the following system of equations applies, again on the space-time domain $Q$,
\begin{equation}
  \label{eq:Finalmodel}
  \begin{split}
      u_t + \nabla \cdot (V^1 \, u) &= -\kappa  u\\
      w_t + \nabla \cdot (V^2 \, w) &= \kappa  u,
  \end{split}
\end{equation}
together with suitable initial and boundary conditions.

\section{Reconstruction of flow parameter fields}\label{sec:reconstruction}

For both models, blood flow parameter fields, i.e., velocities and diffusion or transfer coefficient functions have to be determined from given ultrasound measurements. Here, we use PDE-con\-strained optimization techniques for this task. In this section we discuss the arising optimization problems for both models, derive optimality conditions, and give an algorithm for the reconstruction.

\subsection{Advection-diffusion model}

Denoting the unknown parameter fields by $c = (V, D)$ we consider the Tikhonov-regularized cost functional
\begin{equation}\label{eq:advdiffopt}
  J(u,c) = J(u,(V,D)) = \frac{1}{2} ||u - C_\text{meas}^\delta||_{L_{2}(Q)}^{2} + \frac{\lambda_{1}}{2} ||V||_{L_{2}(\Omega)^d}^{2} + \frac{\lambda_{2}}{2}||D||_{L_{2}(\Omega)}^{2},
\end{equation}
with (potentially noisy) measurements $C_\text{meas}^\delta$.
Assuming local uniqueness and sufficient regularity of the control-to-state mapping, and using the formal Lagrange technique for the reduced cost functional $j(u) = J(u(c),c)$~\cite{troltzsch} we derive the adjoint equation for the adjoint variable $p$
\begin{equation}
  \begin{split}
	    -p_t - \nabla p \cdot V - \nabla \cdot (D \, \nabla p) &= u - C_\text{meas}^\delta \quad \text{ on } Q\\
	    p(\cdot,T) &= 0 \phantom{ - C_\text{meas}^\delta} \quad\ \ \text{ in } \Omega.\\
	  \end{split}
\end{equation}
We assume distributed measurements here due to high frequency and resolution of the imaging modality. Pointwise-in-time measurements could be treated by interpolation, or as jumps in the adjoint equation~\cite{BorziItoKunisch2002}.
For the reduced gradients we obtain
\begin{equation}
  \begin{split}
	    \nabla_V j(u) &= \int_0^T  \lambda_1 \, V + \nabla p \, u \,dt\\
	    \nabla_D j(u) &= \int_0^T  \lambda_2 \, D - \nabla p \cdot \nabla u \,dt,
	  \end{split}
\end{equation}
see Appendix~\ref{sec:adjointAdvecDiff}. 

\subsection{Two-compartment model}

In the two-compartment model the unknown parameters are $c=(V^1, V^2, \kappa)$, and we consider the cost functional
\begin{equation}
  \begin{split}
    J(u,w,c) &= J(u,w,(V^1,V^2,\kappa))\\
    &=  \frac{1}{2} ||u + w - C_\text{meas}^\delta||_{L_{2}(Q)}^{2} 
        + \frac{\lambda_{1}}{2} ||V^1||_{L_{2}(\Omega)^d}^{2} + \frac{\lambda_{2}}{2}||V^2||_{L_{2}(\Omega)^d}^{2}
        + \frac{\lambda_{3}}{2} ||\kappa||_{L_{2}(\Omega)}^{2}.
  \end{split}
\end{equation}
Here, as tracer concentrations in arterial and venous blood cannot be measured separately, only the added concentrations are matched to the measurement data 
 $C_\text{meas}^\delta$. 
Using the control-to-state mapping and assuming sufficient regularity as in the advection-diffusion model, we consider the reduced objective functional $j(u(c),w(c))$ and derive the adjoint equations as
\begin{equation}
  \begin{split}
    -p_t - \nabla p \cdot V^1 + \kappa \, p &= u + w - C_\text{meas}^\delta + \kappa q\\
    -q_t - \nabla q \cdot V^2 &= u + w - C_\text{meas}^{\delta}\\
    p(\cdot,T) =  q(\cdot,T) = 0.
  \end{split}
\end{equation}
Finally, the reduced gradients are given as
\begin{equation}
	\begin{split}
		\nabla_{V} j(u,w) &= \begin{pmatrix}
		    \int_0^T \lambda_1 \, V^1 + \nabla p \, u \,dt\\
			\int_0^T \lambda_2 \, V^2 + \nabla q \, w \,dt
		\end{pmatrix}\\
		\nabla_\kappa j(u,w) &= \int_0^T \lambda_{3}\, \kappa + (q - p) \, u \,dt,
	\end{split}
\end{equation}
see Appendix \ref{sec:adjoint2C}.

\subsection{Implementation}

In order to solve the optimization problems, state and adjoint PDEs have to be approximated numerically.
This is done using a cell-centred uniform grid in space and time. The flux ($\nabla \cdot (V \, u)$) is computed 
via a fifth order weighted essentially non-oscillatory (WENO) scheme \cite{ShuWENONotes}. In time, a implicit solver is needed for the advection-diffusion model
to ensure stability for the diffusive part, thus, a third order IMEX solver is used \cite{IMEXSolver}, combining an explicit Runge-Kutta (ERK)
with a diagonally implicit Runge-Kutta (DIRK) method, where both are strong stability preserving (SSP) and A-stable. The adjoint equation uses the same solver, after the standard time-transformation $\tau = T - t$.
After solving the forward and backward problem, the gradient needs to be calculated in order to update the control variables.
Occurring derivatives and integrals are calculated via a second order central finite differences scheme and quadrature rule, 
respectively. 
We use a splitting approach to update the parameter fields, as adapting the diffusivity fields takes less steps then adapting the velocity.
The optimization starts by adapting the velocity field $V$, before adapting the diffusion $D$. 

For the two-compartment model, the only modification is the time stepping algorithm. 
Because we no longer have a diffusive part, an explicit, total variation diminishing (TVD) third order Runge-Kutta method \cite{ShuWENONotes} is chosen. 
As the parameters have different influence on the cost functional and 
are of different magnitude, we again use splitting, and start by calculating the optimal velocity fields for a fixed transfer coefficient, then adapt $\kappa$, and iterate until optimal parameters are found.
As a step size algorithm the Armijo-algorithm \cite{Armijo} with strong Wolfe conditions \cite{Wolfe} is implemented and the descent 
direction is typically chosen as conjugate gradient method with Dai Yuan parameter \cite{conjugateGradientConvergence} but can also be 
changed to the steepest descent or BFGS method \cite{BFGSNotes}.
We stop the algorithm, when either the gradient norm or parameter updates are sufficiently small. The approach is summarized 
in Algorithm \ref{alg:ParamIdentTwoComp} for the two-compartment model. Modifications for the advection-diffusion model are straightforward. A prototype implementation in 2D is available at \cite{Code}.

\begin{algorithm}[!htp]
  \caption{Determine velocities $V^1$, $V^2$ and transfer coefficient $\kappa$ from  given measured contrast agent concentration $C_\text{meas}^{\delta}$.}
  \label{alg:ParamIdentTwoComp}
  \begin{algorithmic}[1]
    \Require{$V^1_{0}(x)$,$V^2_{0}(x)$, $\kappa_0(x)$, $C_\text{meas}^{\delta}(x,t)$}
    \Ensure{$V^1$, $V^2$, $\kappa$}
    \While{Velocities or $\kappa$ is still adapted}  \Comment{$\ell$-loop}
    	\State $V^{1,(0)} = V^1_{\ell}, V^{2,(0)} = V^2_\ell$
      \While{$||\nabla j_{V} ||> \text{ tol}_1$ and $|J^{(i)}-J^{(i-1)}| > \text{ tol}_2$ and $\alpha_V^{(i)} > \text{ tol}_{\alpha}$} \Comment{$i$-loop}
        \State {Compute $u^{(i)}$, $w^{(i)}$ with given $V^{1,(i)}$, $V^{2,(i)}$ and $\kappa_{\ell}$, via the state equation }
        \State Evaluate objective functional $J^{(i)}$
        \State {Compute $p^{(i)}$, $q^{(i)}$ with given $V^{1,(i)}$, $V^{2,(i)}$, $\kappa_{\ell}$, $u^{(i)}$ and $w^{(i)}$ via the adjoint equation}
        \State {Compute $\nabla j_V$ for velocities}
        \State Find suitable step size via line search: $\alpha_{V}^{(i)}$
        \State Update: $\begin{pmatrix} V^{1,(i+1)}\\ V^{2,(i+1)}\end{pmatrix} = \begin{pmatrix} V^{1,(i)}\\ V^{2,(i)} \end{pmatrix} - \alpha_V^{(i)} \, \nabla j_V $.
      \EndWhile
      \State Optimal velocity fields $V^1_{\ell+1} = V^{1,(i+1)}$ and $V^2_\ell = V^{2,(i+1)}$ for round $\ell$ found. Set $\kappa^{(0)} = \kappa_\ell$
      \While{$||\nabla j_{\kappa} ||> \text{ tol}_1$ and $|J^{(i)}-J^{(i-1)}| > \text{ tol}_2$ and $\alpha_{\kappa}^{(i)} > \text{ tol}_{\alpha}$} \Comment{$i$-loop}
        \State {Compute $u^{(i)}$, $w^{(i)}$ with given $V^{1}_\ell$, $V^{2}_\ell$ and $\kappa^{(i)}$, via the state equation }
       \State Evaluate objective functional $J^{(i)}$
        \State {Compute $p^{(i)}$, $q^{(i)}$ with given $V^{1}_\ell$, $V^{2}_\ell$, $\kappa^{(i)}$, $u^{(i)}$ and $w^{(i)}$ via the adjoint equation}
        \State {Calculate $\nabla j^{\kappa}$ for transfer coefficient}
        \State Find suitable step size via line search: $\alpha_{\kappa}^{(i)}$
        \State Update: $\kappa^{(i+1)} = \kappa^{(i)} - \alpha_{\kappa}^{(i)} \, \nabla_\kappa j$
      \EndWhile
      \State Optimal transfer coefficient $\kappa_\ell= \kappa^{(i+1)} $ for round $\ell$ found
    \EndWhile
    \State Optimal velocity fields $V^1 = V^1_\ell$, $V^2 = V^2_\ell$ and $\kappa = \kappa_\ell$ found
  \end{algorithmic}
\end{algorithm}

\section{Numerical results}\label{sec:results}

We evaluate the reconstruction algorithm on synthetic as well as measurement data. First, we apply the adjoint approach for the advection-diffusion model and compare to an earlier parameter identification approach for this model. To analyze the performance for the two-compartment model, we reconstruct velocities and transfer coefficient in several settings with synthetic data. Finally, using dynamic contrast-enhanced 
ultrasound measurements we compare reconstructions using the advection-diffusion and the two-compartment model. 
We restrict the discussion to 2D in space. While an extension to 3D is in principle straightforward, an optimized 
implementation would be required to avoid prohibitively long computation times and memory requirements, and is left for future work.

\subsection{Adjoint-based parameter estimation for the advection-diffusion model}

For a first evaluation of the adjoint-based method, we compare to an earlier approach by Hristov et al.~\cite{bloodflow}. There, the inverse problem for the advection-diffusion model was formulated as a regularized least-squares problem by differentiating smoothed measurement data and fitting velocity and diffusivity values for each discrete pixel (2D) or voxel (3D). To be more precise, the discrete measurement data was convolved with Gaussian kernels numerically approximating spatial and temporal differentiation, and derivatives of the diffusivity were approximated using backward finite differences. For the model problem, it was assumed that velocities are divergence-free. The resulting linear least-squares problem was solved by a conjugate gradient method. 

We compare the adjoint optimization approach to this using similar synthetic data, albeit in two space dimensions only. Four different combinations of velocity and diffusivity were used to generate synthetic data (GT in Table~\ref{table:benchmarkingStartingValues}). Errors were measured as mean relative errors (MRE),

\begin{equation}\label{mre}
	MRE = \frac{1}{2}\left(\frac{||V - V_{GT}||^2_{L_2(Q)^2}}{||V_{GT}||^2_{L_2(Q)^2}} + \frac{||D - D_{GT}||^2_{L_2(Q)}}{||D_{GT}||^2_{L_2(Q)}}\right).
\end{equation}

For the adjoint approach, the more general model PDE~\eqref{eq:advdiff} was used, not enforcing divergence-free velocity fields. If required, a Leray-type projection~\cite{Chorin1968} could be used for velocity updates to ensure divergence-free velocities, although at higher computational cost.
Steepest descent with strong Wolfe line search was used to solve the optimization problem~\eqref{eq:advdiffopt}. The regularization parameters were set to $\lambda_1 =10^{-4}, \lambda_2 = 10^{-5}$, the
minimum step size to $\text{tol}_{\alpha} = 5\cdot 10^{-7}$, and the tolerances for ending the optimization set to $\text{tol}_1 = \text{tol}_2 = 10^{-5}$.
For discretization,  $90$ time-steps, and a spatial grid with $n_x = n_y = 45$ cells in $x$- and $y$-direction was used. Ground truth (GT) and parameters as well as individual MREs are reported in Table~\ref{table:benchmarkingStartingValues}.

\begin{table}[!ht]
	\centering
	\begin{tabular}{ l c  c  c  c  c  c  c} 
		\toprule
		case & $V = (V_x \  V_y)$  & $D$ & $V^{GT} = (V^{GT}_x \  V^{GT}_y)$ & $D^{GT}$ & MRE \\ 
		\midrule
		1 & $\begin{pmatrix} 0 & 0.5 \end{pmatrix}$ & $0$ & $\begin{pmatrix} 0 & 1 \end{pmatrix}$ & $0.05$ & $0.57$ \\  
		2 & $\begin{pmatrix} 0 & 0.5 \end{pmatrix}$ & $0$ & $\begin{pmatrix} 0 & 1 \end{pmatrix}$ & $0.15$ & $0.52$ \\   
		3 & $\begin{pmatrix} 0 & 1.3 \end{pmatrix}$ & $0$ & $\begin{pmatrix} 0 & 2 \end{pmatrix}$ & $0.05$ & $0.53$ \\  
		4 & $\begin{pmatrix} 0 & 1.3 \end{pmatrix}$ & $0$ & $\begin{pmatrix} 0 & 2 \end{pmatrix}$ & $0.15$ & $0.65$ \\   
		\bottomrule
	\end{tabular}
	\caption{Initial guesses, ground truth, and MRE for comparison to~\cite{bloodflow}.}  
	\label{table:benchmarkingStartingValues}
\end{table}

As expected, the problem with high velocities and high diffusion is the hardest to reconstruct, but even with this straightforward
basic implementation the adjoint approach gives a combined MRE of $0.5675$. This is indication of improvement compared to the  least-squares approach~\cite{bloodflow}, where a combined total MRE of $0.681$ was reported, and shows feasibility of the adjoint approach.Advantages of the continuous adjoint approach include not only the efficient and flexible computation of gradients, e.g., allowing for different discretization schemes for state and adjoint equations, but also allows to de-couple the discretization of the PDE from the resolutions of the imaging modality. This is especially beneficial for unequally sampled or sparse data points, as well as moving regions of interest as occuring in DCE-US~\cite{Tiyarattanachai2022}.

\subsection{Two-compartment model: synthetic data}

To investigate performance of the adjoint approach for the two-compartment model, we generate synthetic data on the space-time domain $Q = \Omega \times (0,T) = [1, 3] \times [1, 3] \times (0,1)$. We discretize 
with $40 \times 40 $ cells in space, and use 120 time steps. We pose homogeneous Neumann boundary conditions, assuming that 
the contrast bolus, modeled by a Gaussian shaped initial condition for the arterial compartment, is contained fully in the 
domain at initial time. Boundary conditions are prescribed using three ghost cells on every side of the domain. The parameters 
with which the data sets have been created numerically using the two-compartment model
can be found in Table \ref{table:SyntheticData}. In order to model the 
behavior of blood flow in an organ, the arterial blood velocity is fast at the beginning and zero towards the end of the domain, 
while the venous blood velocity is zero at the beginning and faster towards the end of the domain. 
Furthermore, the velocities in $y$-direction are chosen to model a spreading out in the arterial concentrations and a contracting
in the venous concentrations around the middle of the domain. To analyze the influence of the conversion domain on the propagation of contrast agent, the transfer coefficient $\kappa$ is set to a constant value on a wide domain, that is either prescribed (WTD) or needs to be identified (WTD-s), as well as on a narrow domain (NTD), and $\kappa = 0$ outside these domains.  
Furthermore, we test the reconstruction on problem WTD with additive noise, to demonstrate robustness of the algorithm.
For this, we add $5\%$ and $10\%$ Gaussian noise \cite{gaussianNoise}, with respect to the maximum concentration, to the
artificial data, using $\mu = 0$ for both cases and $\sigma = 0.15$ and $\sigma = 0.3$ for the $5$ and $10\%$ noise, respectively.

\begin{table}[!ht]
    \centering
    \begin{tabular}{ p{2cm} p{5cm} p{5cm} } 
        \toprule
        parameters & WTD / WTD-s & NTD\\ 
        \midrule
        $V^1_x$  & $\begin{cases} -3x + 7, & x\leq 2.3 \\ \phantom{-}0, & \text{else} \end{cases} $ & $ \begin{cases} -3x + 7, & x\leq 2.1 \\ \phantom{-}0, & \text{else} \end{cases} $ \\
        $V^1_y$  & $\begin{cases} \phantom{-}0.3, &  y < 2 \text{ and } x \in [1.5,2] \\ -0.3, & y \geq 2 \text{ and } x \in [1.5,2] \\ \phantom{-}0, & \text{else} \end{cases}$  & $\begin{cases} \phantom{-}0.3,  & y < 2 \text{ and } x \in [1.5,2] \\ -0.3, & y \geq 2 \text{ and } x \in [1.5,2] \\ \phantom{-}0, & \text{else} \end{cases} $ \\
        \midrule
        $V^2_x$  & $\frac{5}{4}x - \frac{5}{4}$ & $\frac{5}{4}x - \frac{5}{4}$ \\
        $V^2_y$  & $ \begin{cases} \phantom{-}0.3, & y < 2 \text{ and } x \in [2,2.5] \\ -0.3, & y \geq 2 \text{ and } x \in [2,2.5] \\ \phantom{-}0, & \text{else} \end{cases} $ & $ \begin{cases} \phantom{-}0.3, & y < 2 \text{ and } x \in [2,2.5] \\ -0.3, & y \geq 2 \text{ and } x \in [2,2.5] \\ \phantom{-}0, & \text{else} \end{cases} $ \\
        \midrule
        $\kappa$ & $7$ & $9$ \\
        $\kappa$-domain & $[1.5, 2.5] \times [1, 3]$ & $[1.8, 2.2] \times [1, 3]$ \\
        \midrule
        $u(x,y, 0)$ & $3 \exp\left(-\frac{(x-1.1)^2 + (y-2)^2}{0.02}\right)$ & $3 \exp\left(-\frac{(x-1.1)^2 + (y-2)^2}{0.02}\right)$ \\
        $w(x,y,0)$ & $0$ & $0$ \\
        \bottomrule
    \end{tabular}
    \caption{Parameters with which the synthetic data has been created.}
    \label{table:SyntheticData}
\end{table}

\begin{table}[!ht]
    \centering
    \begin{tabular}{ l  c   c   c   c   c   c  c } 
        \toprule
        name & $V^1 = (V^1_x \  V^1_y)$ & $V^2 = (V^2_x \  V^2_y)$ & $\kappa$ &  $\{\lambda_i, i=1,2,3\}$ & $\text{tol}_1$ & $\text{tol}_2$ & $\text{tol}_{\alpha}$\\ 
        \midrule
        WTD/WTD-s & $\begin{pmatrix} 2 & 0 \end{pmatrix}$ & $\begin{pmatrix} 2.2 & 0 \end{pmatrix}$ & $18$ & $\{10^{-4},10^{-4},10^{-5}\}$ & $10^{-5}$ & $10^{-7}$ & $5\cdot 10^{-5}$ \\  
        NTD & $\begin{pmatrix} 1.8 & 0 \end{pmatrix}$ & $\begin{pmatrix} 2 & 0 \end{pmatrix}$ & $16$ & $\{10^{-4},10^{-4},10^{-4}\}$& $10^{-4}$ & $10^{-5}$ & $5\cdot 10^{-5}$ \\  
        \bottomrule
    \end{tabular}
    \caption{Initialization and parameters for the reconstructions of synthetic data.}  
    \label{table:SyntheticDataStarting}
\end{table}

For the synthetic test case, the reconstruction is run with the same resolution as the data to avoid interpolation.
The initial guesses and parameters that were used to run the parameter identifications can be found in Table \ref{table:SyntheticDataStarting}.
Here, $\text{tol}_{1}$ is the stopping criterion that ensures that the gradient norm is small enough, $\text{tol}_{2}$
ensures that the decrease in the cost functional is sufficient and $\text{tol}_{\alpha}$ the parameter
that ensures that the chosen step size is large enough. 
For the problems with added noise the same starting parameters as for problem {WTD} were chosen,
except for problem {WTD $10\%$} where we had to change $\text{tol}_{2}=10^{-6}$ in order to perform sufficiently many optimization iterations.
Furthermore, all results where obtained with the conjugate gradient method as descent direction, the descent direction is
reset to the steepest descent direction if the difference between the previous and current search direction 
$||s_{k}-s_{k-1}||^2_{L_2(Q)} \leq 5\cdot 10^{-4}$. During the strong Wolfe line search we use $c_1 = 0.1$ and $c_2 = 10^{-4}$.

\begin{table}[!ht]
    \centering
    \begin{tabular}{ p{2.5cm} p{2.5cm} p{2.5cm} p{2.7cm} p{2.3cm}} 
        \toprule
        name & $||\nabla j_{v}||^2_{L_2(Q)}$ & $||\nabla j_{\kappa}||^2_{L_2(Q)}$ &  $J_\text{final}\newline  (J_\text{initial})$ & $|\kappa-\kappa_\text{exact}|$ \newline (rel. err.) \\ 
        \midrule
        WTD       & $0.012$ & $9.3\cdot 10^{-4}$ & $0.0082$ ($0.0738$) & $0.24$ ($3.4\%$)\\
        WTD 5\%  & $0.0032$ & $4\cdot 10^{-4}$ & $0.0126$ ($0.0769$) & $0.04$ ($0.57\%$) \\
        WTD 10\% & $29.43$ & $1.7\cdot 10^{-4}$ & $0.0289$ ($0.0837$) & $0.49$ ($7\%$) \\
         WTD-s   & $0.00057$ & $1 \cdot 10^{-7}$ & $0.0127$ ($0.06$) & -\\
        \midrule
        NTD       & $0.0003$ & $0.0017$ & $0.004$ ($0.1438$) & $0.14$ ($1.6\%$)\\
        \bottomrule
    \end{tabular}
    \caption{Main indicators of success for reconstructions of the synthetic test for the two-compartment model.}  
    \label{table:resultNoise}
\end{table}

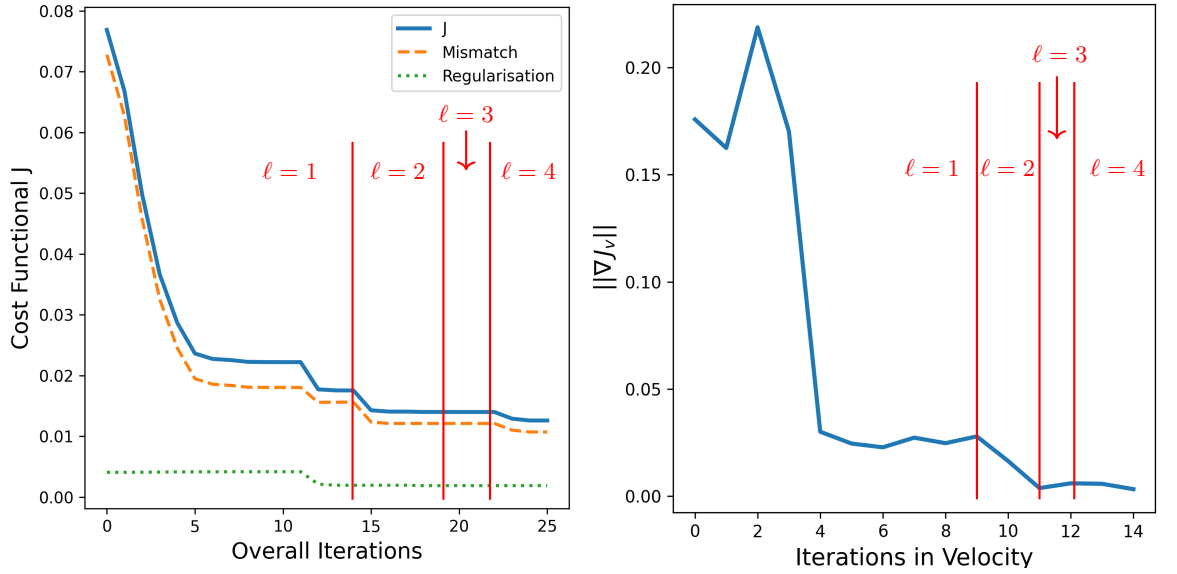
\begin{figure}[!htp]
	\centering
	\begin{subfigure}[b]{0.49\textwidth}
		\centering
		\begin{tikzpicture}
			
			\node[anchor=south west, inner sep=0] (img) 
			{\includegraphics[width=\textwidth]{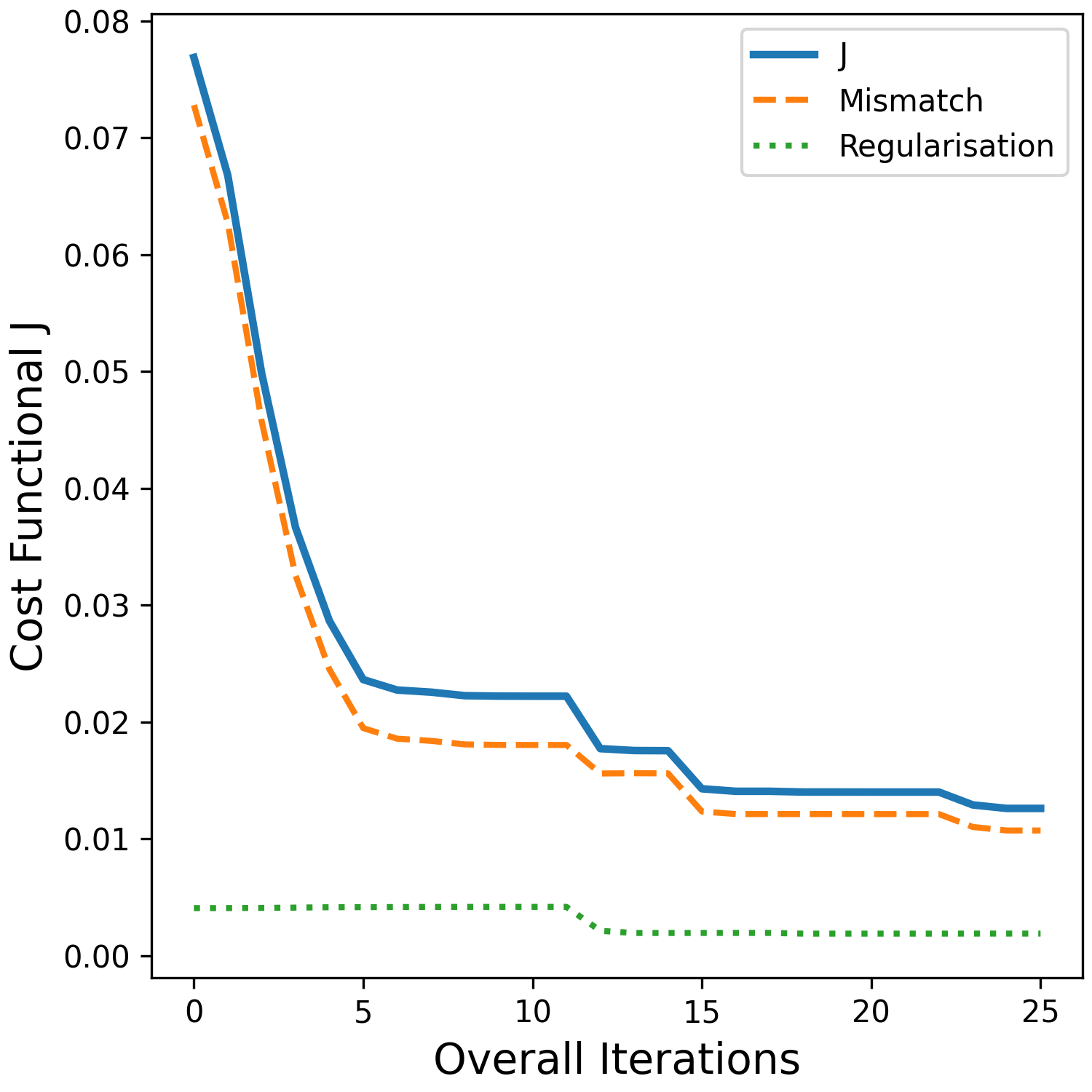}};
			
			\begin{scope}[x={(img.south east)}, y={(img.north west)}]
				\draw[red, thick] (0.61,0.12) -- (0.61,0.75);
				\draw[red, thick] (0.77,0.12) -- (0.77,0.75);
				\draw[red, thick] (0.852,0.12) -- (0.852,0.75);
				\node[font=\bfseries\small, text=red]
				at (0.5,0.7) {$\ell = 1$};
				
				\node[font=\bfseries\small, text=red]
				at (0.69,0.7) {$\ell = 2$};
				
				\node[font=\bfseries\small, text=red]
				at (0.81,0.8) {$\ell = 3$};
				\draw[->, red, thick] (0.81,0.77) -- (0.81,0.7);
				
				\node[font=\bfseries\small, text=red]
				at (0.92,0.7) {$\ell = 4$};
				
			\end{scope}
		\end{tikzpicture}
		\label{fig:cF_noise5}
	\end{subfigure}
	\hfill
	\begin{subfigure}[b]{0.5\textwidth}
		\centering
		\begin{tikzpicture}
			
			\node[anchor=south west, inner sep=0] (img) 
			{\includegraphics[width=\textwidth]{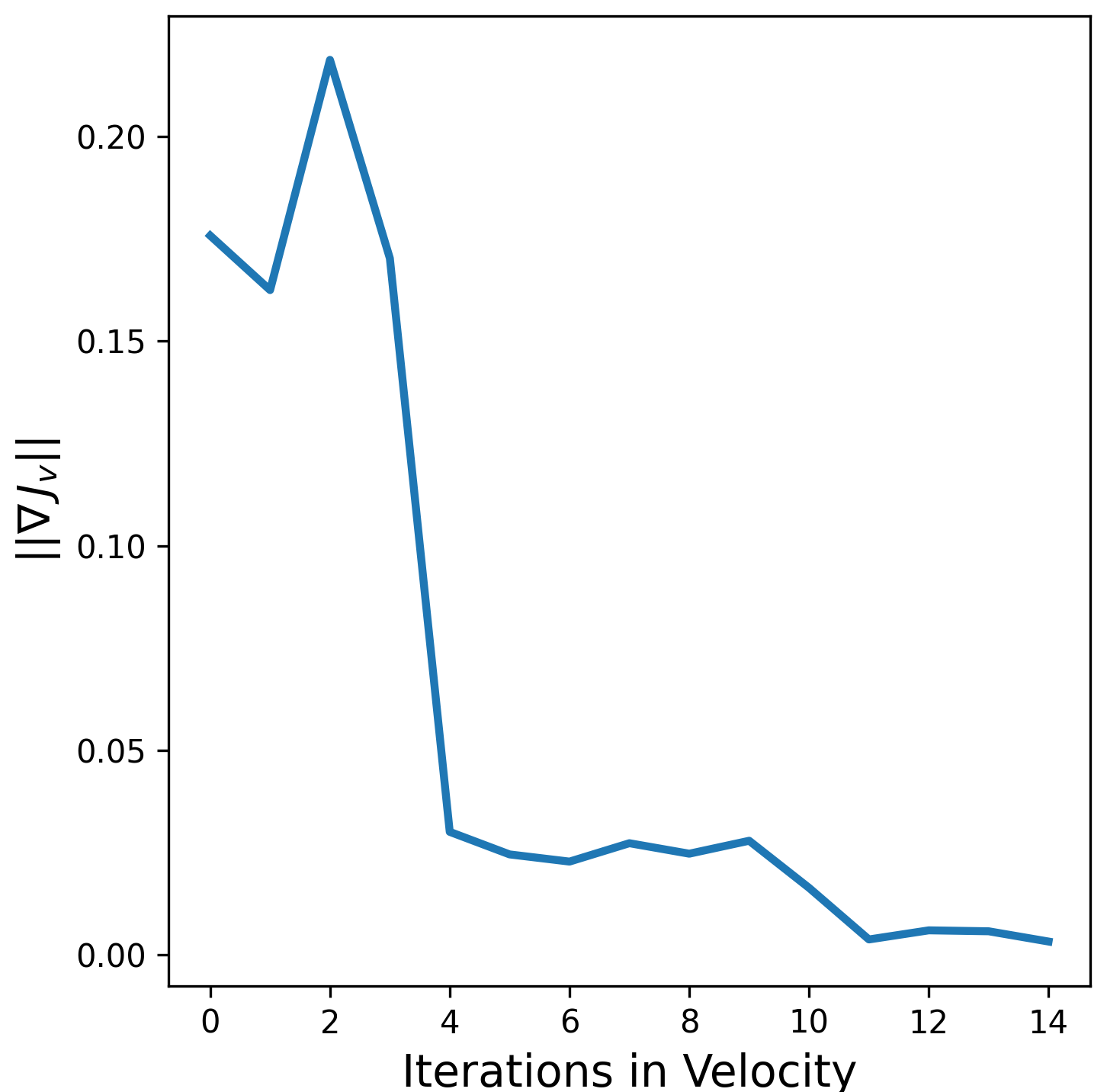}};
			
			\begin{scope}[x={(img.south east)}, y={(img.north west)}]
				\draw[red, thick] (0.677,0.12) -- (0.677,0.85);
				\draw[red, thick] (0.785,0.12) -- (0.785,0.85);
				\draw[red, thick] (0.845,0.12) -- (0.845,0.85);
				\node[font=\bfseries\small, text=red]
				at (0.6,0.7) {$\ell = 1$};
				
				\node[font=\bfseries\small, text=red]
				at (0.73,0.7) {$\ell = 2$};
				
				\node[font=\bfseries\small, text=red]
				at (0.82,0.9) {$\ell = 3$};
				\draw[->, red, thick] (0.815,0.86) -- (0.815,0.75);
				
				\node[font=\bfseries\small, text=red]
				at (0.92,0.7) {$\ell = 4$};
				
			\end{scope}
		\end{tikzpicture}
		\label{fig:grad_vel_noise5}
	\end{subfigure}
	\caption{Cost functional and its contributions from the data mismatch and regularization term (left) as well as gradient norm of velocities (right; only showing velocity iterations) for the reconstruction of problem WTD with $5 \%$ noise and $4$ optimization rounds, showing the progress of the optimization.}
	\label{fig:IndicNoise5}
\end{figure}

In all cases we can see a convergence in cost functional value, which is exemplarily plotted for {WTD $5\%$}
in Figure \ref{fig:IndicNoise5}, together with the progress of the norm of the gradient with respect to the velocity over the velocity iterations.

The final cost functional values for all problems can be found in Table \ref{table:resultNoise}.
The final cost functional values are lower for the problems without noise. As noise was 
added on the whole domain, the mismatch between reconstructed concentrations and the artificial data is higher, because it also
appears outside the area of the tracer transport, i.e., in areas where the concentrations $u+w = 0$ and $C_\text{meas}^\delta$ is only noise. A remedy for this would be the use of multiplicative noise, e.g., speckle noise~\cite{Michailovich2003,Goodman2020}.

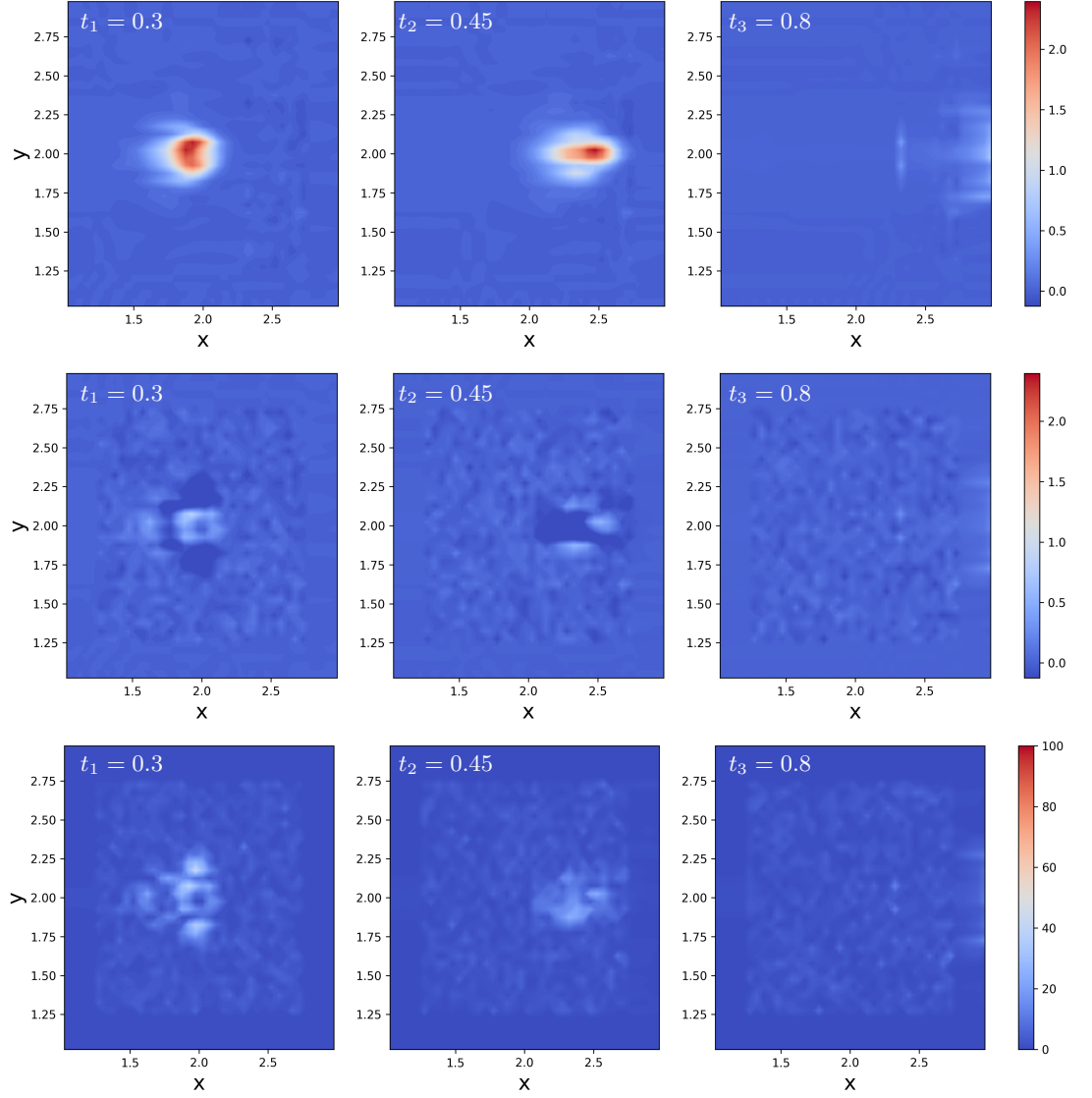
\begin{figure}[!ht]
	\centering
	\begin{subfigure}[c]{\textwidth}
		\centering
		\begin{tikzpicture}
			
			\node[anchor=south west, inner sep=0] (img) 
			{\includegraphics[width=0.95\textwidth]{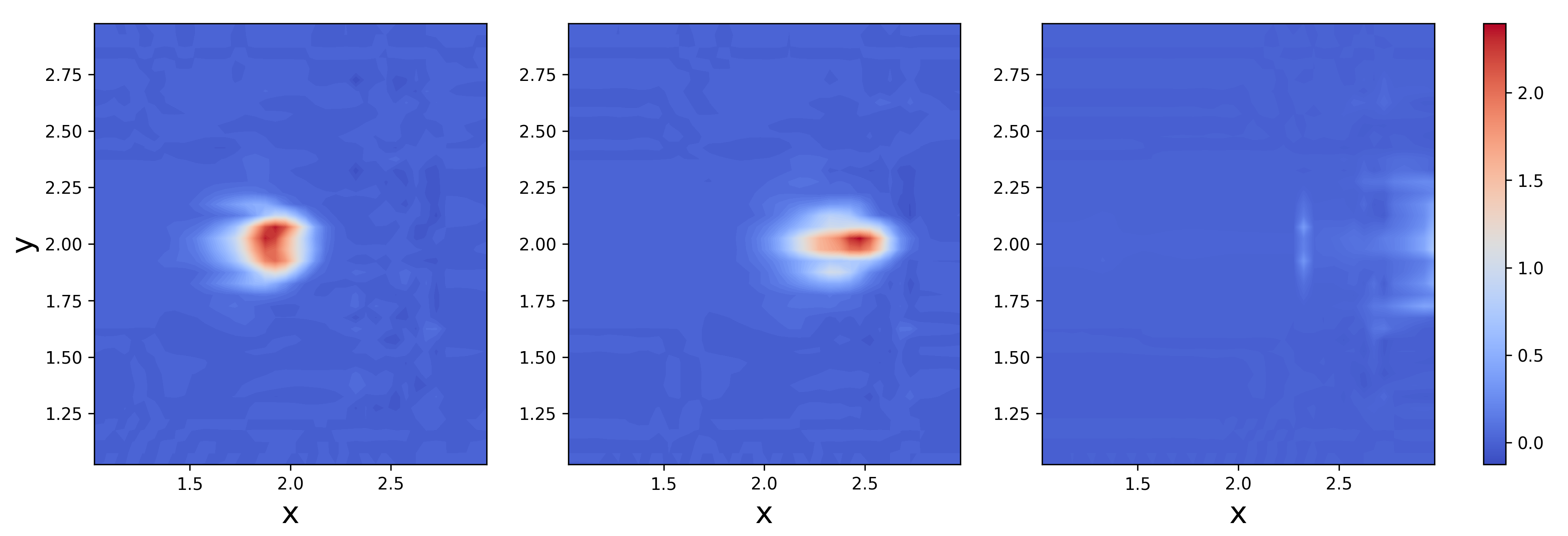}};
			
			\begin{scope}[x={(img.south east)}, y={(img.north west)}]
				
				\node[font=\bfseries\small, text=white]
				at (0.11,0.9) {$t_1 = 0.3$};
				
				\node[font=\bfseries\small, text=white]
				at (0.41,0.9) {$t_2 = 0.45$};
				
				\node[font=\bfseries\small, text=white]
				at (0.71,0.9) {$t_3 = 0.8$};
				
			\end{scope}
		\end{tikzpicture}
	\end{subfigure}
	\begin{subfigure}[c]{\textwidth}
		\centering
		\begin{tikzpicture}
			
			\node[anchor=south west, inner sep=0] (img) 
			{\includegraphics[width=0.95\textwidth]{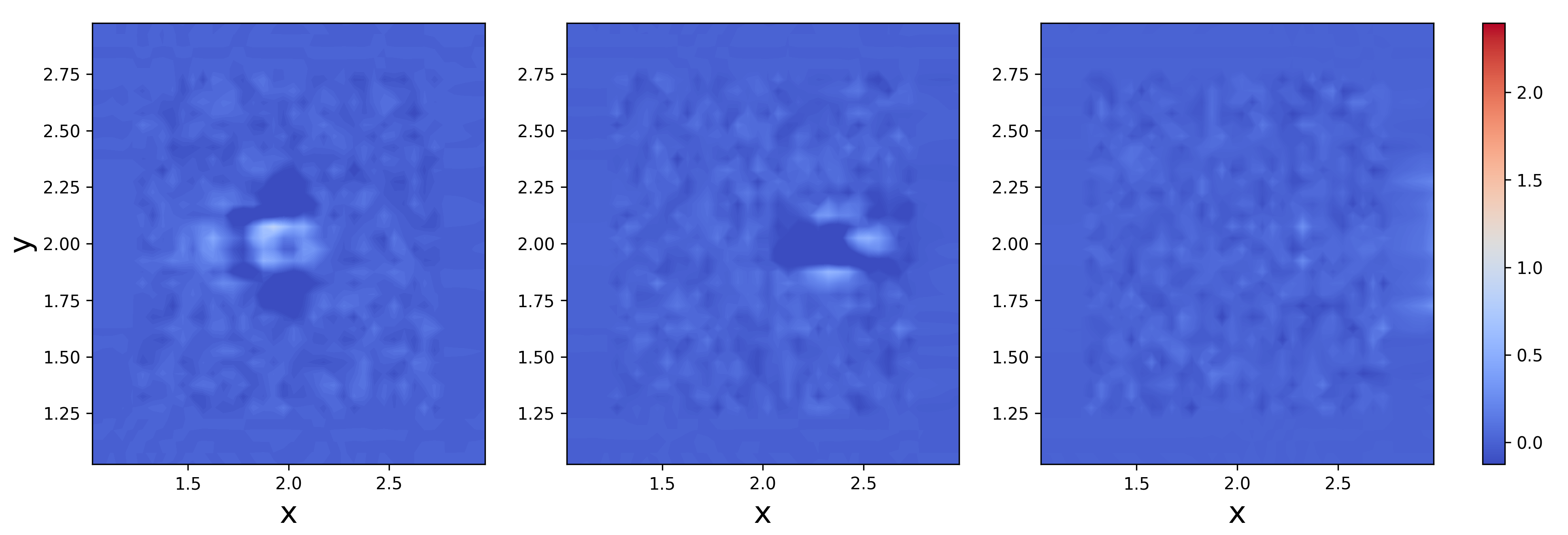}};
			
			\begin{scope}[x={(img.south east)}, y={(img.north west)}]
				
				\node[font=\bfseries\small, text=white]
				at (0.11,0.9) {$t_1 = 0.3$};
				
				\node[font=\bfseries\small, text=white]
				at (0.41,0.9) {$t_2 = 0.45$};
				
				\node[font=\bfseries\small, text=white]
				at (0.71,0.9) {$t_3 = 0.8$};
				
			\end{scope}
		\end{tikzpicture}
	\end{subfigure}
	\begin{subfigure}[c]{\textwidth}
		\centering
		\begin{tikzpicture}
			
			\node[anchor=south west, inner sep=0] (img) 
			{\includegraphics[width=0.95\textwidth]{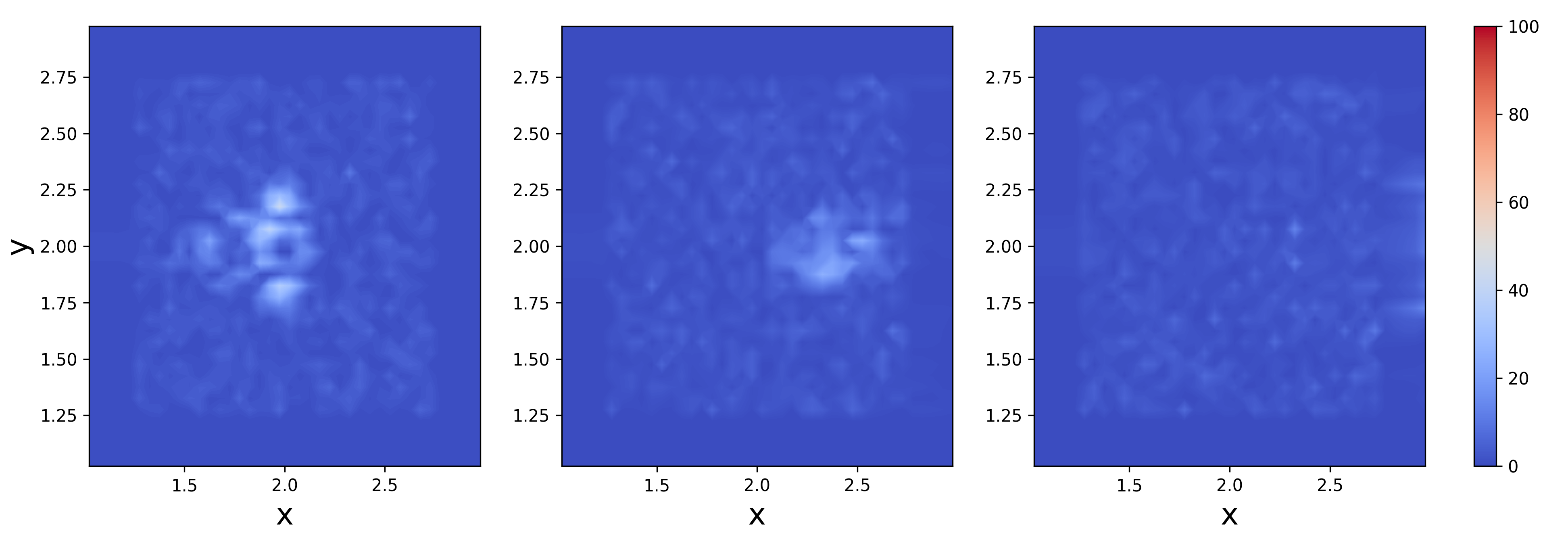}};
			
			\begin{scope}[x={(img.south east)}, y={(img.north west)}]
				
				\node[font=\bfseries\small, text=white]
				at (0.11,0.9) {$t_1 = 0.3$};
				
				\node[font=\bfseries\small, text=white]
				at (0.41,0.9) {$t_2 = 0.45$};
				
				\node[font=\bfseries\small, text=white]
				at (0.71,0.9) {$t_3 = 0.8$};
				
			\end{scope}
		\end{tikzpicture}
	\end{subfigure}
	\caption{Reconstructed concentrations (\textit{top}), the mismatch of reconstructed concentrations and artificial data (\textit{middle}) and the relative error in percent (\textit{bottom}) of problem WTD with $5 \%$ noise after $4$ optimization rounds plotted at three time steps.}
	\label{fig:concRecNoise5}
\end{figure}

\begin{figure}[!htp]
	\centering
	\begin{subfigure}[c]{\textwidth}
		\centering
		\begin{tikzpicture}
			
			\node[anchor=south west, inner sep=0] (img) 
			{\includegraphics[width=0.95\textwidth]{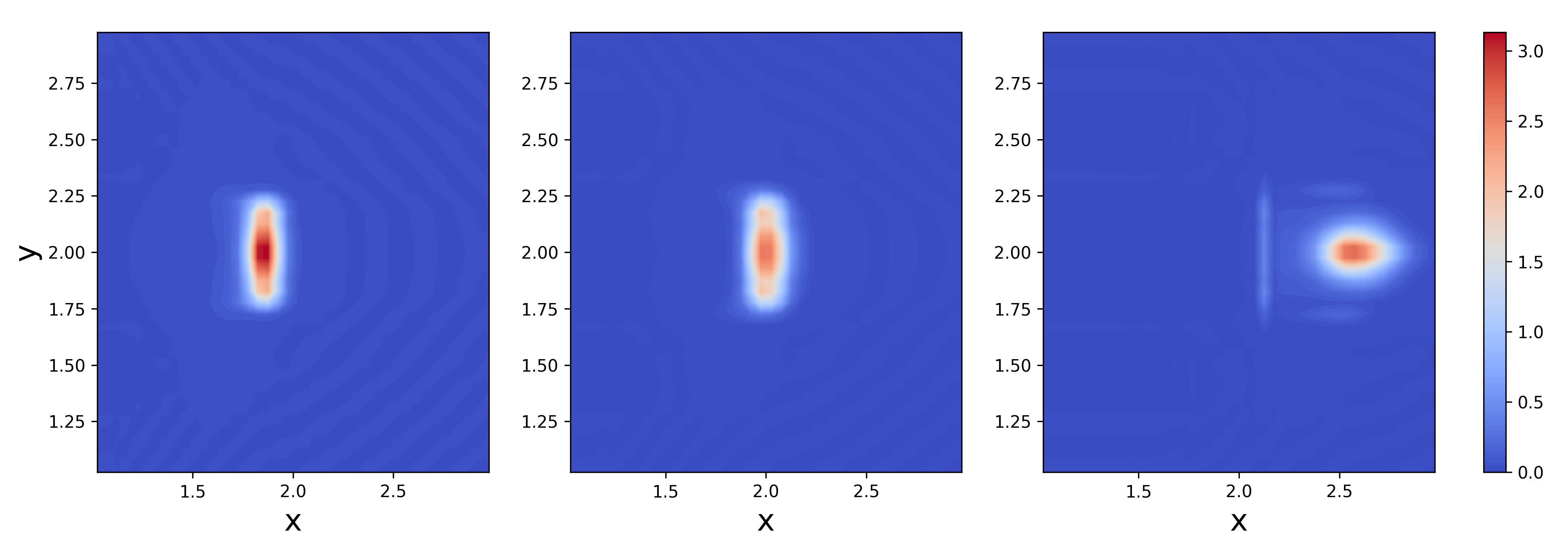}};
			
			\begin{scope}[x={(img.south east)}, y={(img.north west)}]
				
				\node[font=\bfseries\small, text=white]
				at (0.11,0.9) {$t_1 = 0.3$};
				
				\node[font=\bfseries\small, text=white]
				at (0.41,0.9) {$t_2 = 0.4$};
				
				\node[font=\bfseries\small, text=white]
				at (0.715,0.9) {$t_3 = 0.65$};
				
			\end{scope}
		\end{tikzpicture}
	\end{subfigure}
	\hfill
	\begin{subfigure}[c]{\textwidth}
		\centering
		\begin{tikzpicture}
			
			\node[anchor=south west, inner sep=0] (img) 
			{\includegraphics[width=0.95\textwidth]{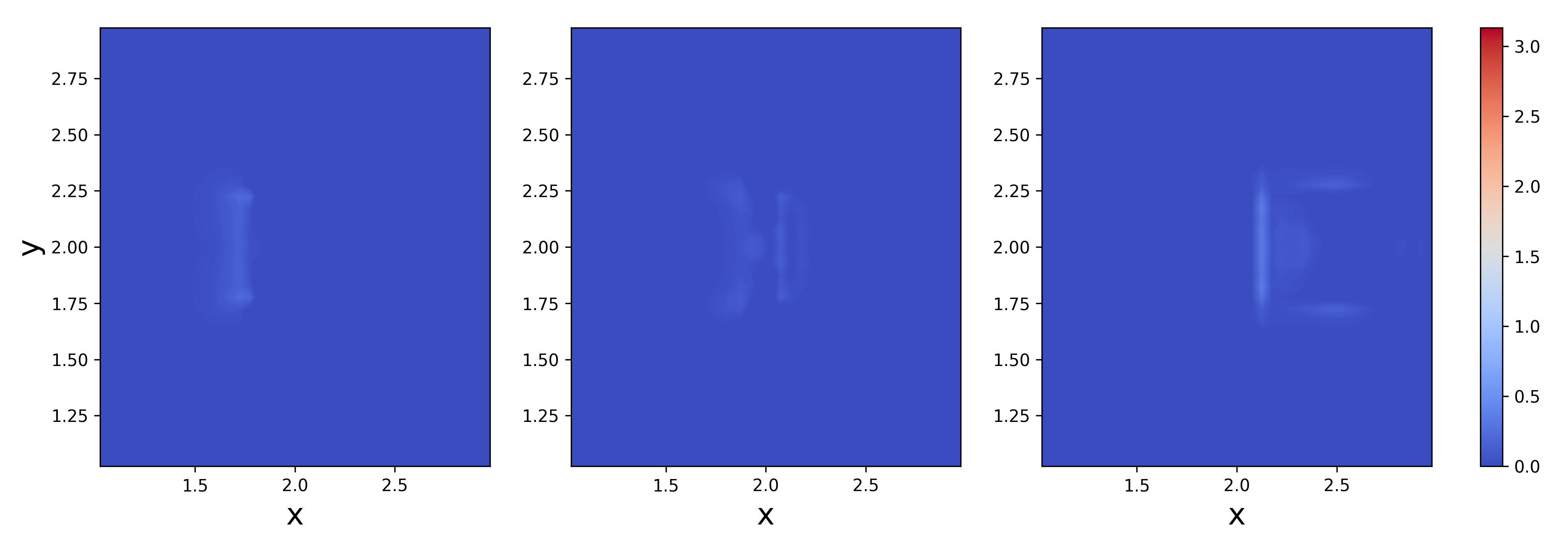}};
			
			\begin{scope}[x={(img.south east)}, y={(img.north west)}]
				
				\node[font=\bfseries\small, text=white]
				at (0.11,0.9) {$t_1 = 0.3$};
				
				\node[font=\bfseries\small, text=white]
				at (0.41,0.9) {$t_2 = 0.4$};
				
				\node[font=\bfseries\small, text=white]
				at (0.715,0.9) {$t_3 = 0.65$};
				
			\end{scope}
		\end{tikzpicture}
	\end{subfigure}
	\begin{subfigure}[c]{\textwidth}
		\centering
		\begin{tikzpicture}
			
			\node[anchor=south west, inner sep=0] (img) 
			{\includegraphics[width=0.95\textwidth]{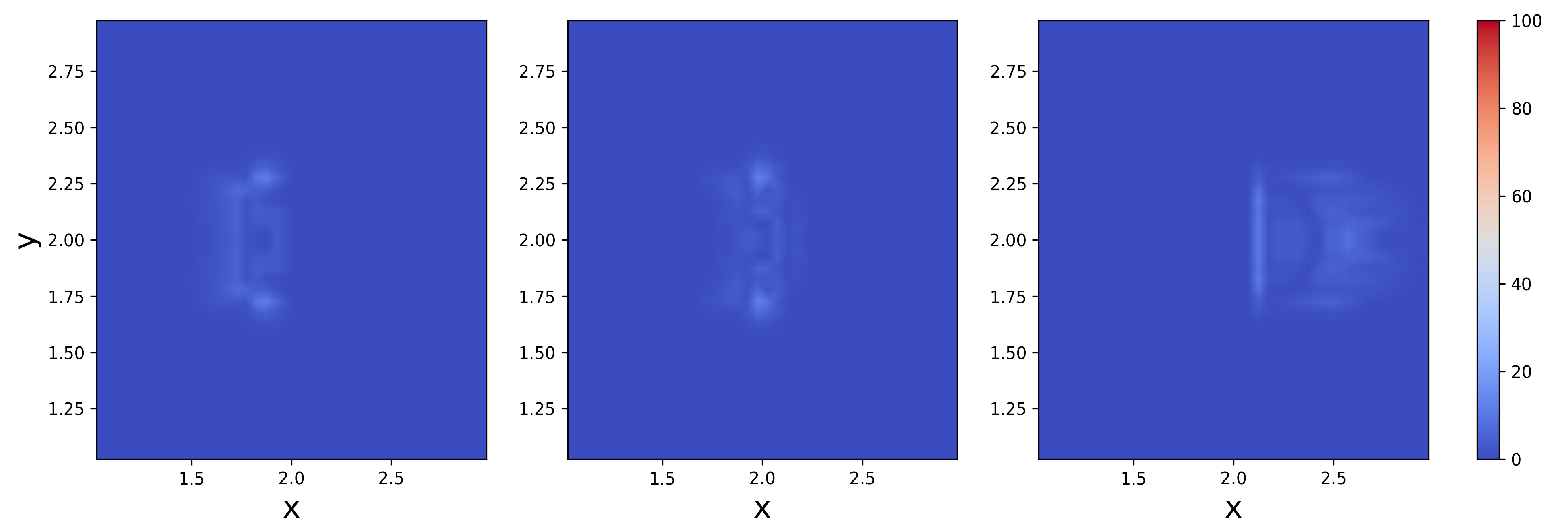}};
			
			\begin{scope}[x={(img.south east)}, y={(img.north west)}]
				
				\node[font=\bfseries\small, text=white]
				at (0.11,0.9) {$t_1 = 0.3$};
				
				\node[font=\bfseries\small, text=white]
				at (0.41,0.9) {$t_2 = 0.4$};
				
				\node[font=\bfseries\small, text=white]
				at (0.715,0.9) {$t_3 = 0.65$};
				
			\end{scope}
		\end{tikzpicture}
	\end{subfigure}
	\caption{Reconstructed concentrations (\textit{top}), mismatch between reconstructed concentrations and artificial data (\textit{middle}) and the relative error in percent (\textit{bottom}) of problem NTD after $24$ optimization rounds.}
	\label{fig:SpreadOutSimReconstruction2}
\end{figure}
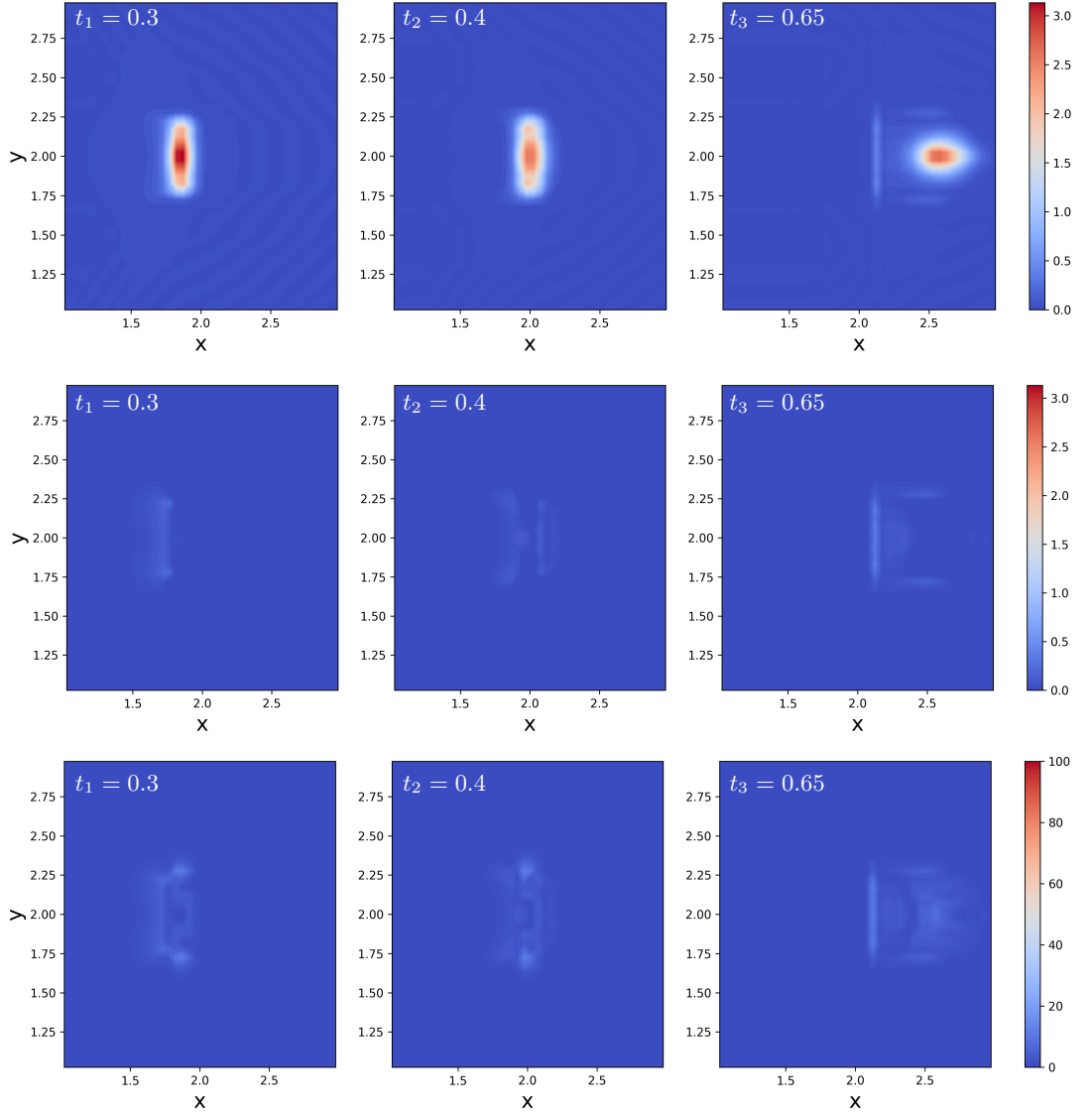

\begin{figure}[!htp]
	\centering
	\begin{tikzpicture}
		\node[anchor=south west, inner sep=0] (img) 
		{\includegraphics[width=0.95\textwidth]{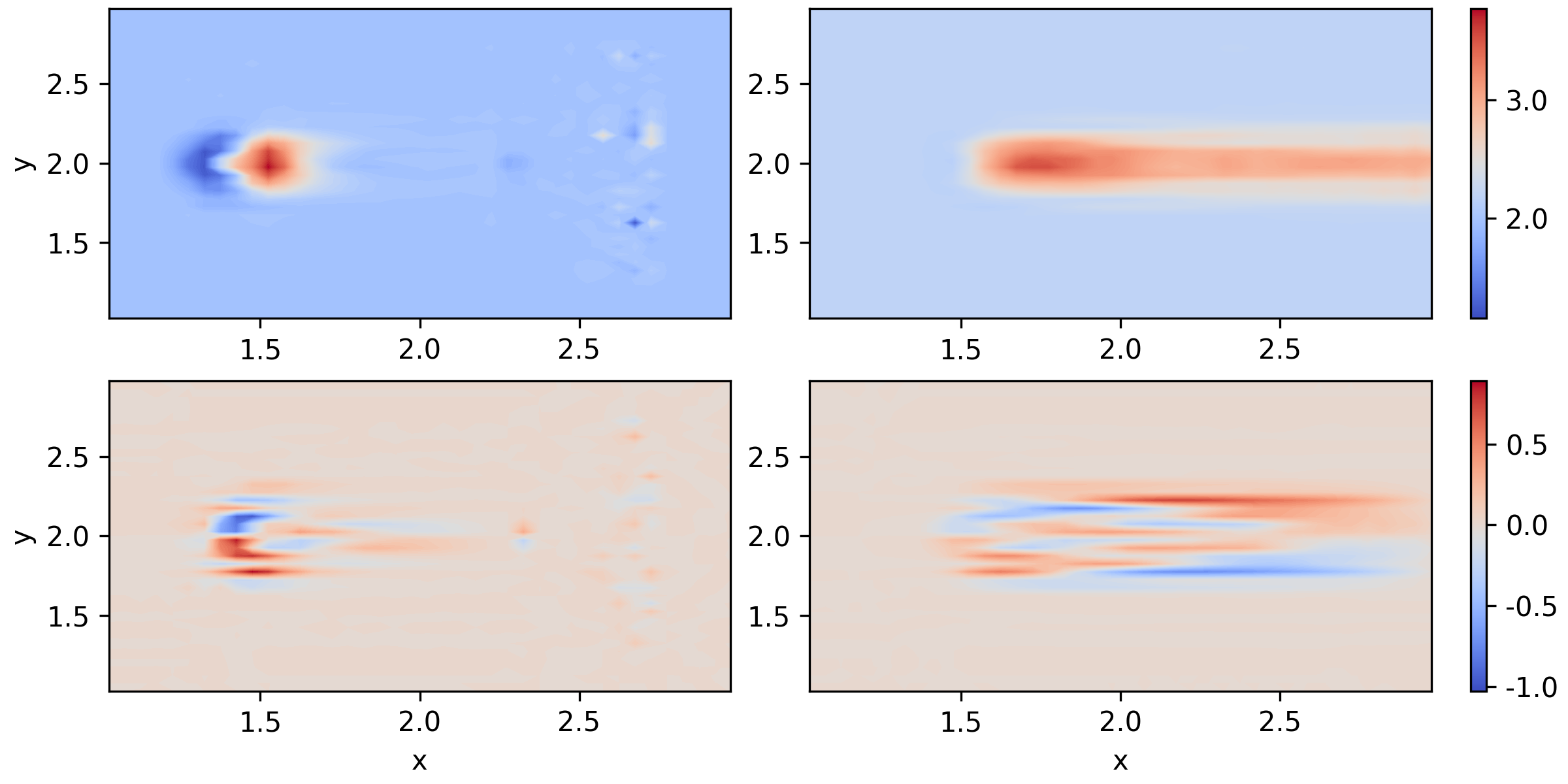}};
		
		\begin{scope}[x={(img.south east)}, y={(img.north west)}]
			
			\node[font=\bfseries, text=black]
			at (0.1,0.94) {$V^1_x$};
			
			\node[font=\bfseries, text=black]
			at (0.545,0.94) {$V^2_x$};
			
			\node[font=\bfseries, text=black]
			at (0.1,0.465) {$V^1_y$};
			
			\node[font=\bfseries, text=black]
			at (0.545,0.465) {$V^2_y$};
			
		\end{scope}
	\end{tikzpicture}
	\caption{Reconstructed velocities of problem WTD with $5 \%$ noise and $4$ optimization rounds.}
	\label{fig:vel_noise5}
\end{figure}
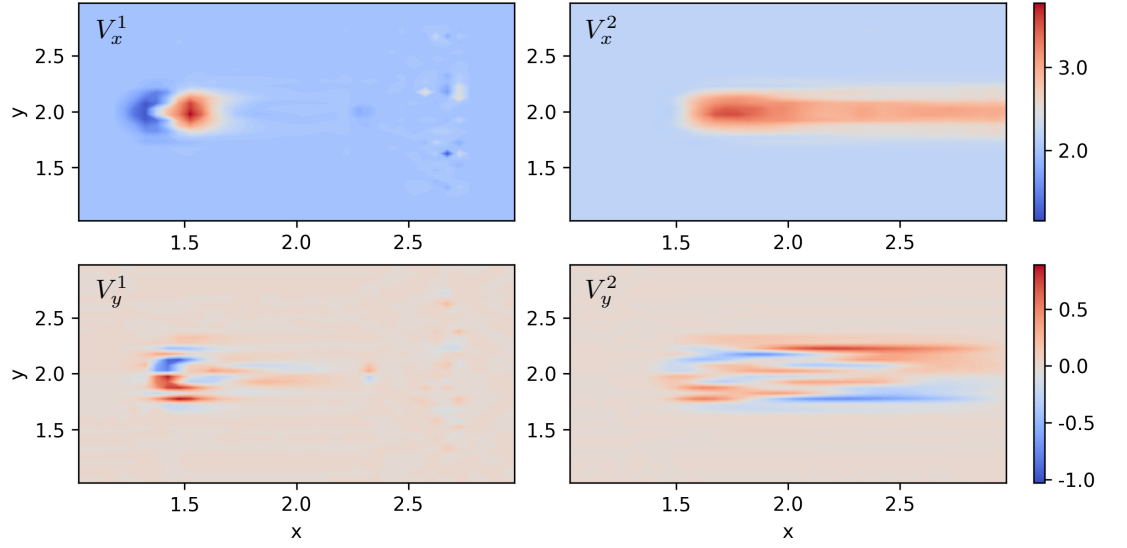

\begin{figure}[!htp]
	\centering
	\begin{tikzpicture}
		\node[anchor=south west, inner sep=0] (img) 
		{\includegraphics[width=0.95\textwidth]{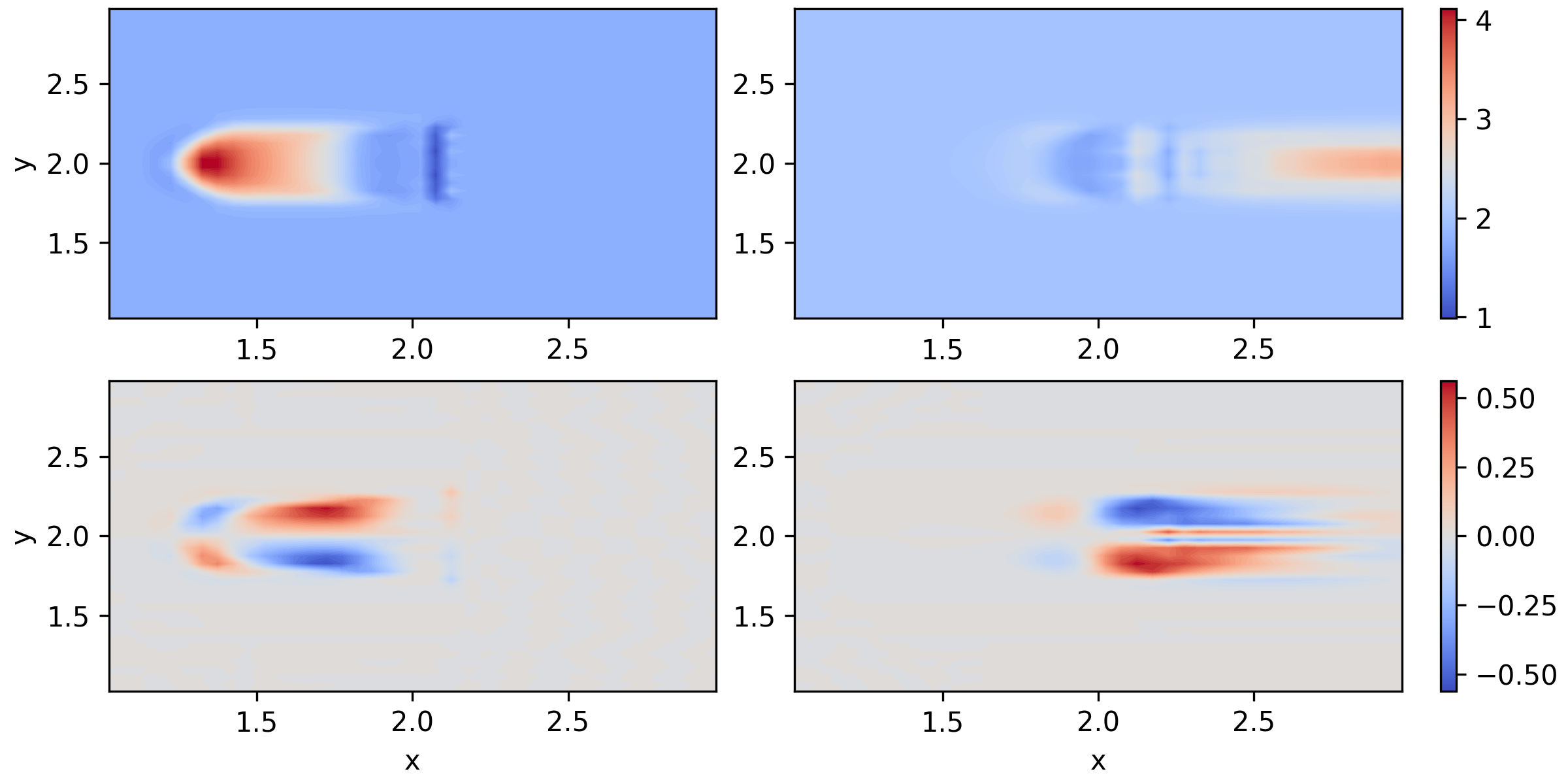}};
		
		\begin{scope}[x={(img.south east)}, y={(img.north west)}]
			
			\node[font=\bfseries, text=black]
			at (0.1,0.94) {$V^1_x$};
			
			\node[font=\bfseries, text=black]
			at (0.535,0.94) {$V^2_x$};
			
			\node[font=\bfseries, text=black]
			at (0.1,0.465) {$V^1_y$};
			
			\node[font=\bfseries, text=black]
			at (0.535,0.465) {$V^2_y$};
			
		\end{scope}
	\end{tikzpicture}
	\caption{Reconstructed velocities of problem NTD after $24$ optimization rounds.}
	\label{fig:reconstructedVel2}
\end{figure}

Furthermore, the gradient norm for the velocities and transfer coefficient, see Table \ref{table:resultNoise}, have been reduced 
in all cases but problem {WTD $10\%$}, thus indicating that the reconstruction is successful. 
The simulated concentrations using the reconstructed parameters are shown
exemplarily for problem {WTD $5\%$} and {NTD} in Figures \ref{fig:concRecNoise5} and \ref{fig:SpreadOutSimReconstruction2}.
It is visible, that both reconstructions are able to fit the movement as well as the spreading out of the bolus at the beginning of the simulation
and the merging towards the end of the simulation and spatial domain. This indicates that  velocities $V^1, V^2$
and transfer coefficient $\kappa$ value are reconstructed sufficiently, see Figures \ref{fig:vel_noise5} and \ref{fig:reconstructedVel2}. 
Comparing both tracer flows shows how good the reconstruction can deal with
different settings for the conversion domain and successfully reconstructing those. 

In the case of $10\%$ noise, the velocities could not be identified sufficiently, as indicated by the large gradient norm at termination. Here,
the bolus movement is compromised by the noise, spreading out and merging of the bolus is 
barely visible in the reconstruction. Furthermore, the concentrations have a wrong shift towards the  
lower part of the domain. This is the case even though the reconstructed velocities 
look similar to the previous ones. Despite this, the transfer coefficient (on a fixed conversion domain), is quite well identified, and is only $0.49$ off the exact value $\kappa = 7$.

\begin{figure}[!ht]
	\centering
	\begin{subfigure}[c]{\textwidth}
		\centering
		\begin{tikzpicture}
			
			\node[anchor=south west, inner sep=0] (img) 
			{\includegraphics[width=0.95\textwidth]{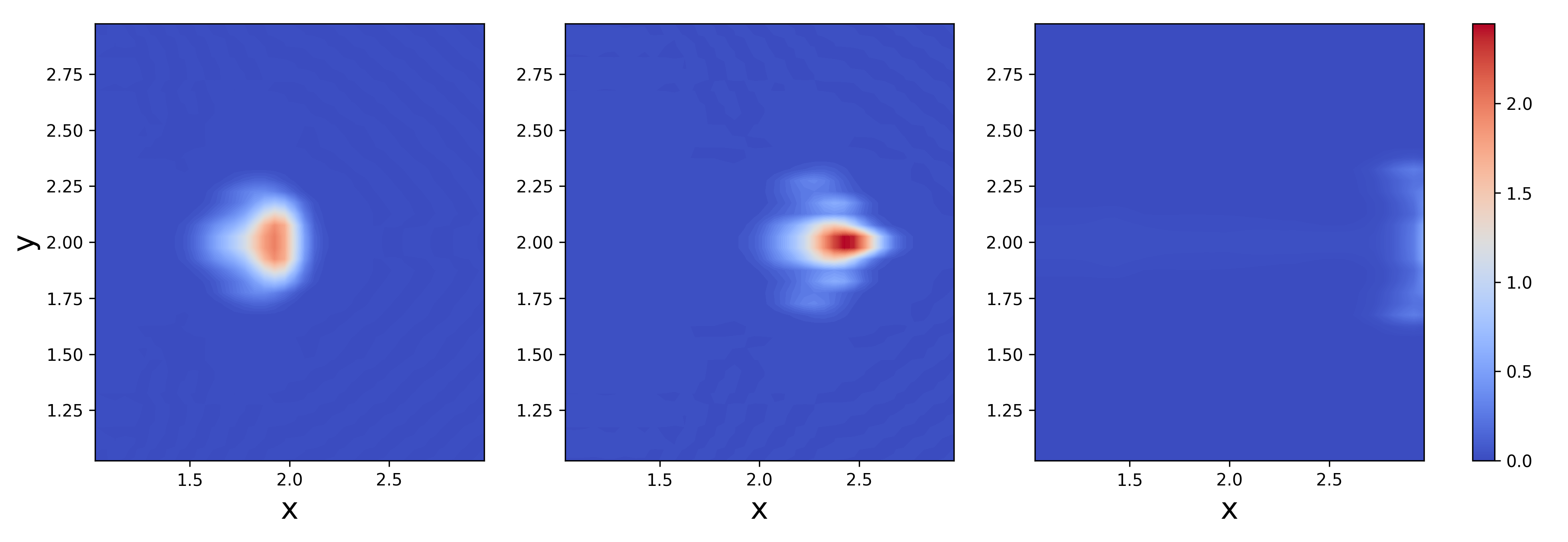}};
			
			\begin{scope}[x={(img.south east)}, y={(img.north west)}]
				
				\node[font=\bfseries\small, text=white]
				at (0.11,0.9) {$t_1 = 0.3$};
				
				\node[font=\bfseries\small, text=white]
				at (0.41,0.9) {$t_2 = 0.45$};
				
				\node[font=\bfseries\small, text=white]
				at (0.71,0.9) {$t_3 = 0.8$};
				
			\end{scope}
		\end{tikzpicture}
		
	\end{subfigure}
	\begin{subfigure}[c]{\textwidth}
		\centering
		\begin{tikzpicture}
			
			\node[anchor=south west, inner sep=0] (img) 
			{\includegraphics[width=0.95\textwidth]{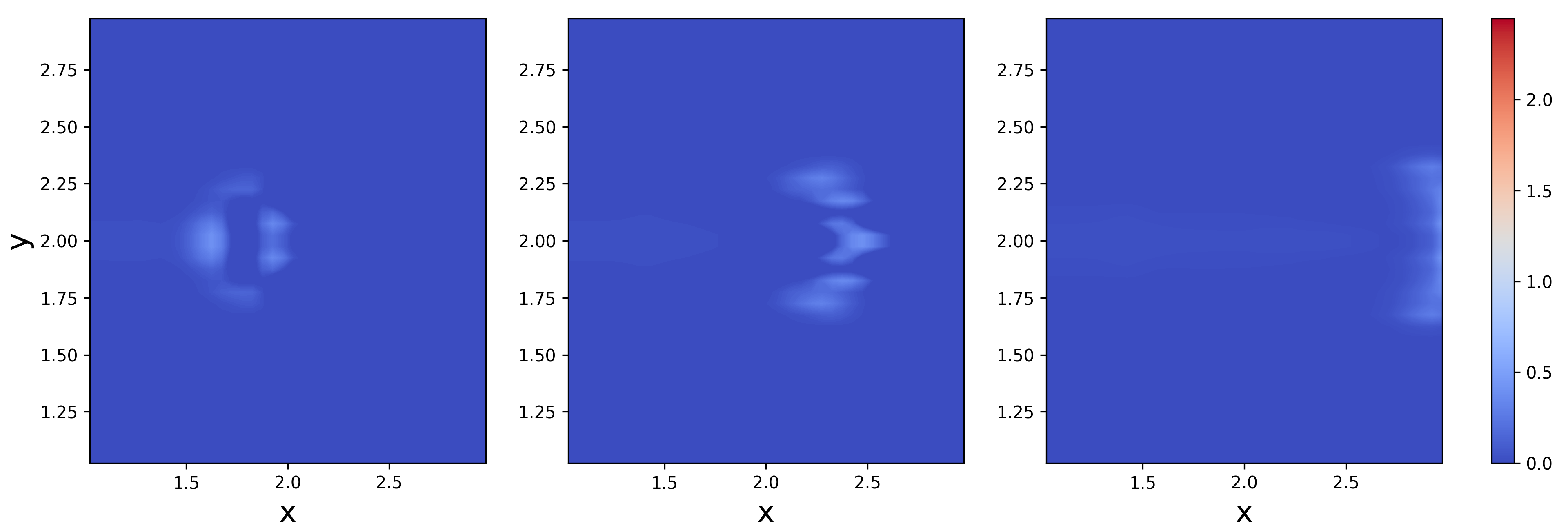}};
			
			\begin{scope}[x={(img.south east)}, y={(img.north west)}]
				
				\node[font=\bfseries\small, text=white]
				at (0.11,0.9) {$t_1 = 0.3$};
				
				\node[font=\bfseries\small, text=white]
				at (0.41,0.9) {$t_2 = 0.45$};
				
				\node[font=\bfseries\small, text=white]
				at (0.71,0.9) {$t_3 = 0.8$};
				
			\end{scope}
		\end{tikzpicture}
	\end{subfigure}
	\begin{subfigure}[c]{\textwidth}
		\centering
		\begin{tikzpicture}
			
			\node[anchor=south west, inner sep=0] (img) 
			{\includegraphics[width=0.95\textwidth]{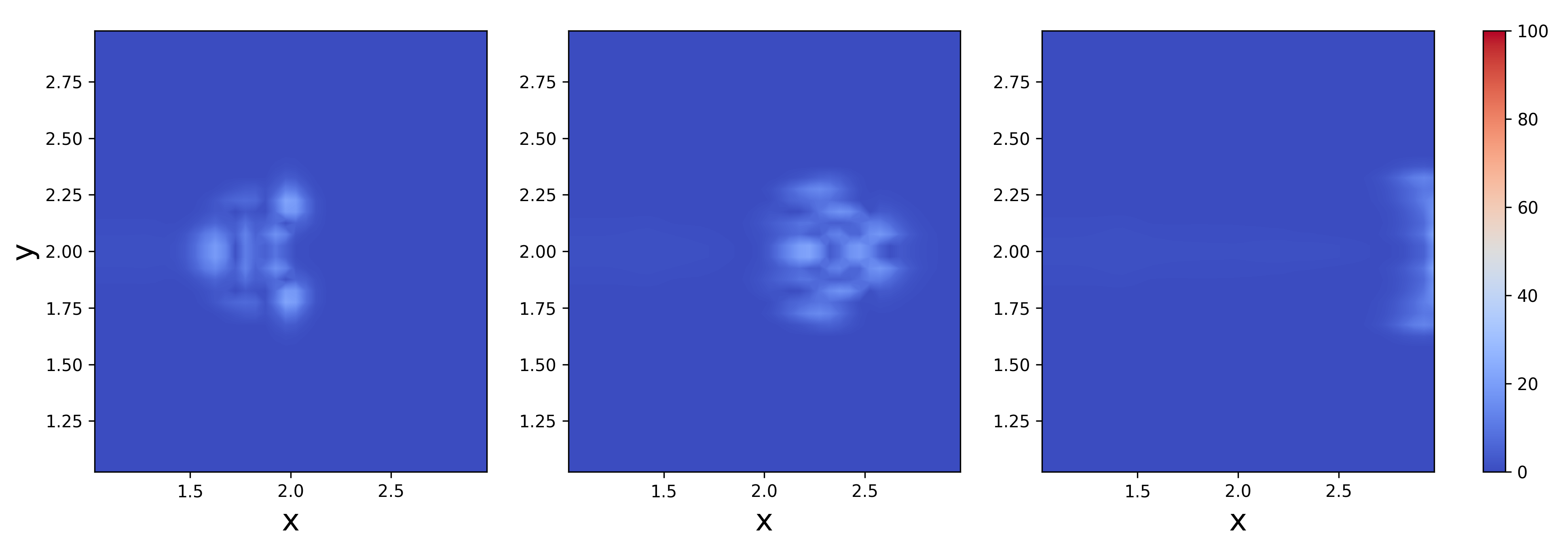}};
			
			\begin{scope}[x={(img.south east)}, y={(img.north west)}]
				
				\node[font=\bfseries\small, text=white]
				at (0.11,0.9) {$t_1 = 0.3$};
				
				\node[font=\bfseries\small, text=white]
				at (0.41,0.9) {$t_2 = 0.45$};
				
				\node[font=\bfseries\small, text=white]
				at (0.71,0.9) {$t_3 = 0.8$};
				
			\end{scope}
		\end{tikzpicture}
		
	\end{subfigure}
	\caption{Reconstructed concentrations (\textit{top}), the mismatch of reconstructed concentrations and artificial data (\textit{middle}) and the relative error in percent (\textit{bottom}) of problem WTD with space-dependent $\kappa$ after $3$ optimization rounds plotted at three time steps.}
	\label{fig:concWTDSpace}
\end{figure}

\begin{figure}[!htp]
	\centering
	\begin{subfigure}{0.74\textwidth}
	\begin{tikzpicture}
		\node[anchor=south west, inner sep=0] (img) 
		{\includegraphics[width=0.95\textwidth]{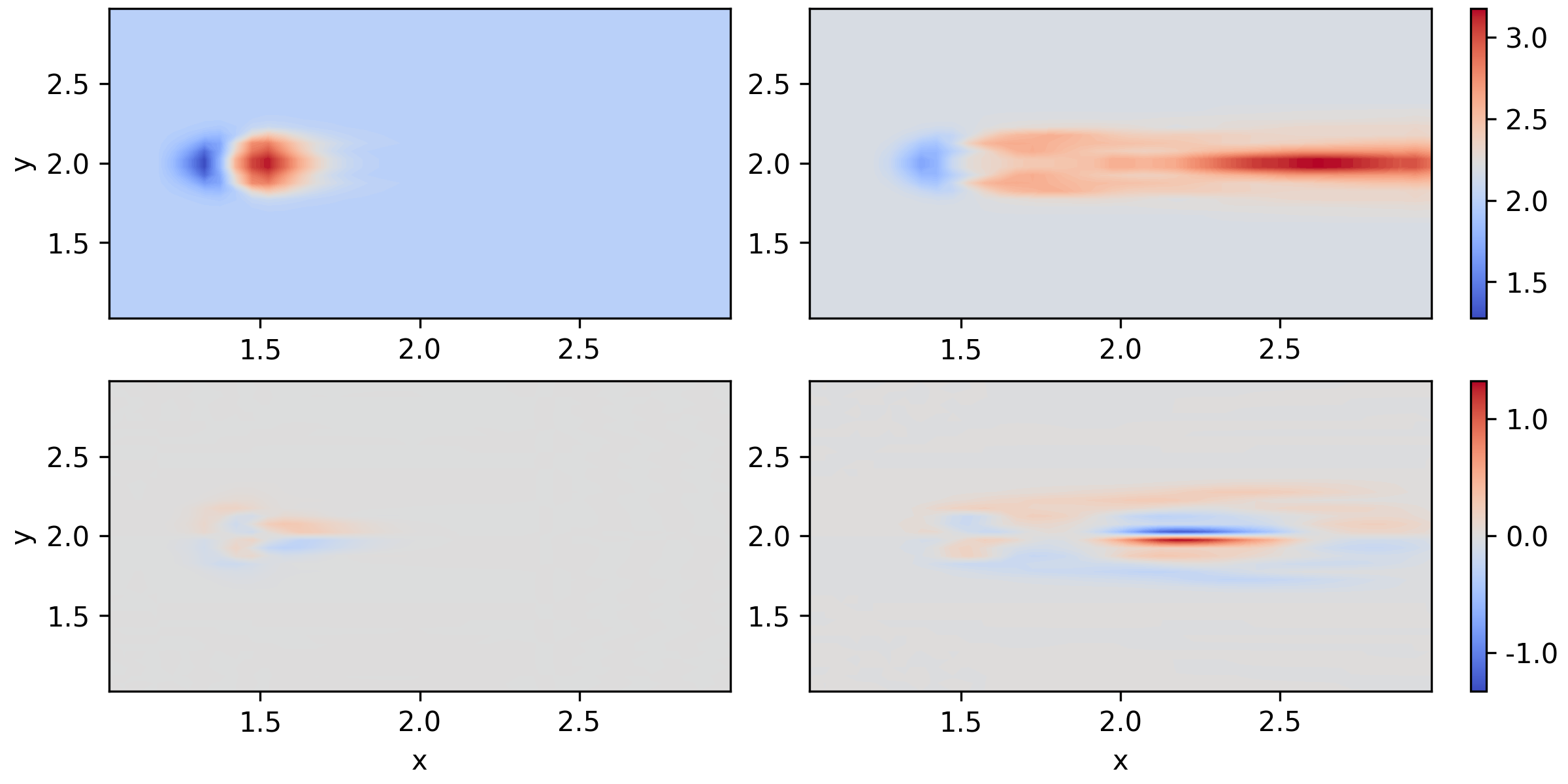}};
		
		\begin{scope}[x={(img.south east)}, y={(img.north west)}]
			
			\node[font=\bfseries, text=black]
			at (0.1,0.94) {$V^1_x$};
			
			\node[font=\bfseries, text=black]
			at (0.545,0.94) {$V^2_x$};
			
			\node[font=\bfseries, text=black]
			at (0.1,0.465) {$V^1_y$};
			
			\node[font=\bfseries, text=black]
			at (0.545,0.465) {$V^2_y$};
			
		\end{scope}
	\end{tikzpicture}
	\end{subfigure}
	\begin{subfigure}{0.25\textwidth}
		\begin{tikzpicture}
		\node[anchor=south west, inner sep=0] (img) 
		{\raisebox{1cm}{\includegraphics[width=\textwidth]{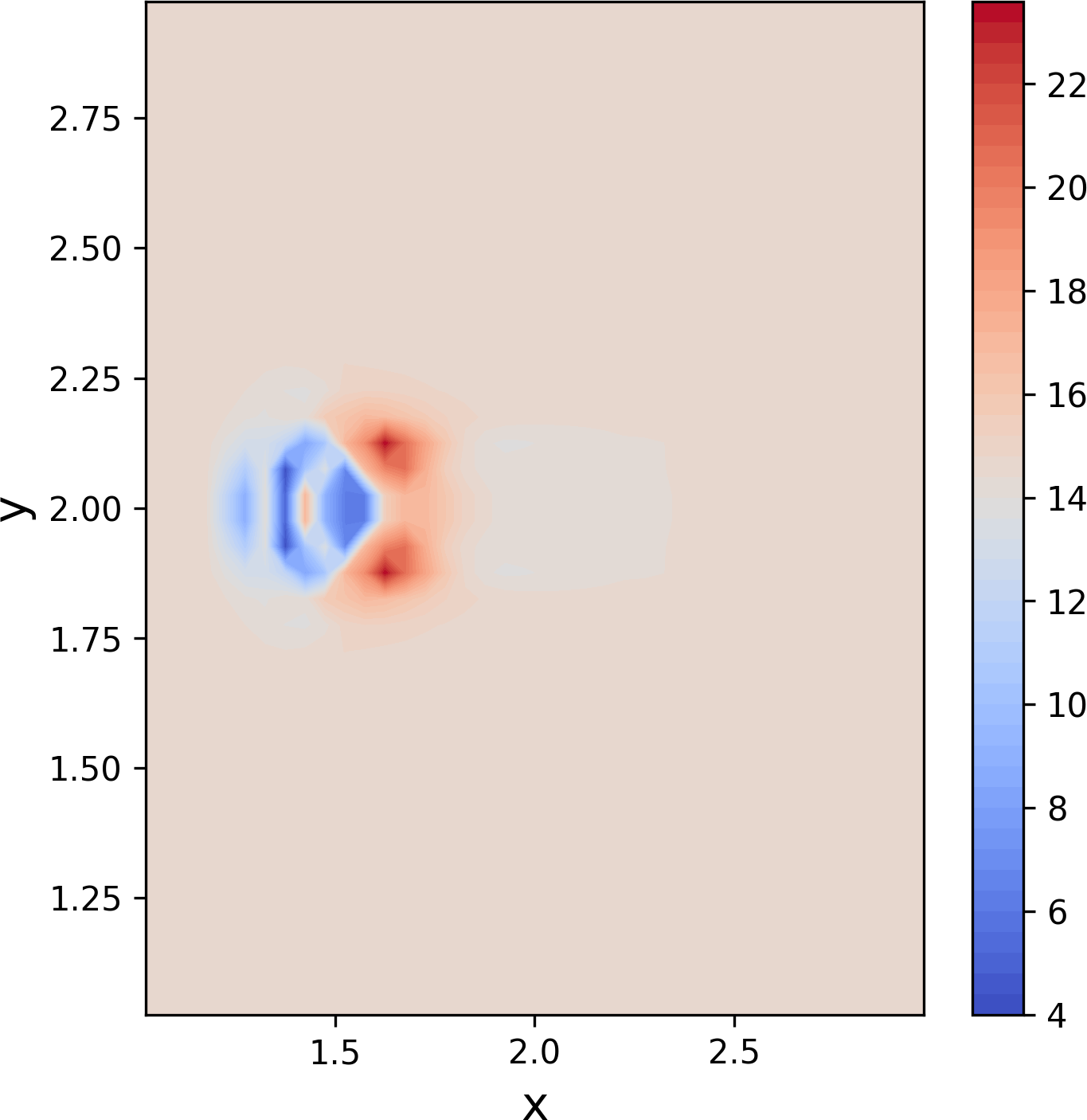}}};
		\begin{scope}[x={(img.south east)}, y={(img.north west)}]
			\node[font=\bfseries, text=black]
			at (0.185,0.96) {$\kappa$};
		\end{scope}
	\end{tikzpicture}
	\end{subfigure}
	
	\caption{Reconstructed velocities and transfer coefficient function of problem WTD-s  after $3$ optimization rounds.}
	\label{fig:vel_WTD_space}
\end{figure}

As a good estimate of the conversion domain is typically not available a priori, , instead of prescribing this domain, we reconstruct a space-dependent transfer coefficient $\kappa(x,y)$ using steepest descent (WTD-s). Despite less prior knowledge being used in this experiment, the reconstruction results are satisfactory, as shown in the concentrations (Figure~\ref{fig:concWTDSpace}) and reconstructed velocities and transfer coefficient (Figure~\ref{fig:vel_WTD_space}). 

To conclude, the parameter identification successfully reconstructs the perfusion parameters for the synthetic test cases, even if moderate noise is introduced into the data. 

\subsection{Measurement data}

We evaluate the reconstruction on
experimental animal data acquired for  previous studies on
3D Dynamic Contrast-Enhanced Ultrasound~\cite{elkaffasSpatialCharacterizationTumor2020,bloodflow}. The data was measured on spatial domain $\Omega = 25.8mm \times 13.2mm \times 23.1mm$ 
with $86$, $44$ and $77$ voxels in $x$-, $y$- and $z$-direction, respectively. The time is given in seconds 
and taken on time domain $T = [0, 26.789s]$ with $50$ time steps, resulting in a frame rate of a little under $2$ Hertz.

\subsubsection{Data preprocessing.}\label{sec:dataProcess}

\begin{figure}
    \centering
    \begin{subfigure}[b]{0.5\textwidth}
        \centering
        \includegraphics[width=\textwidth]{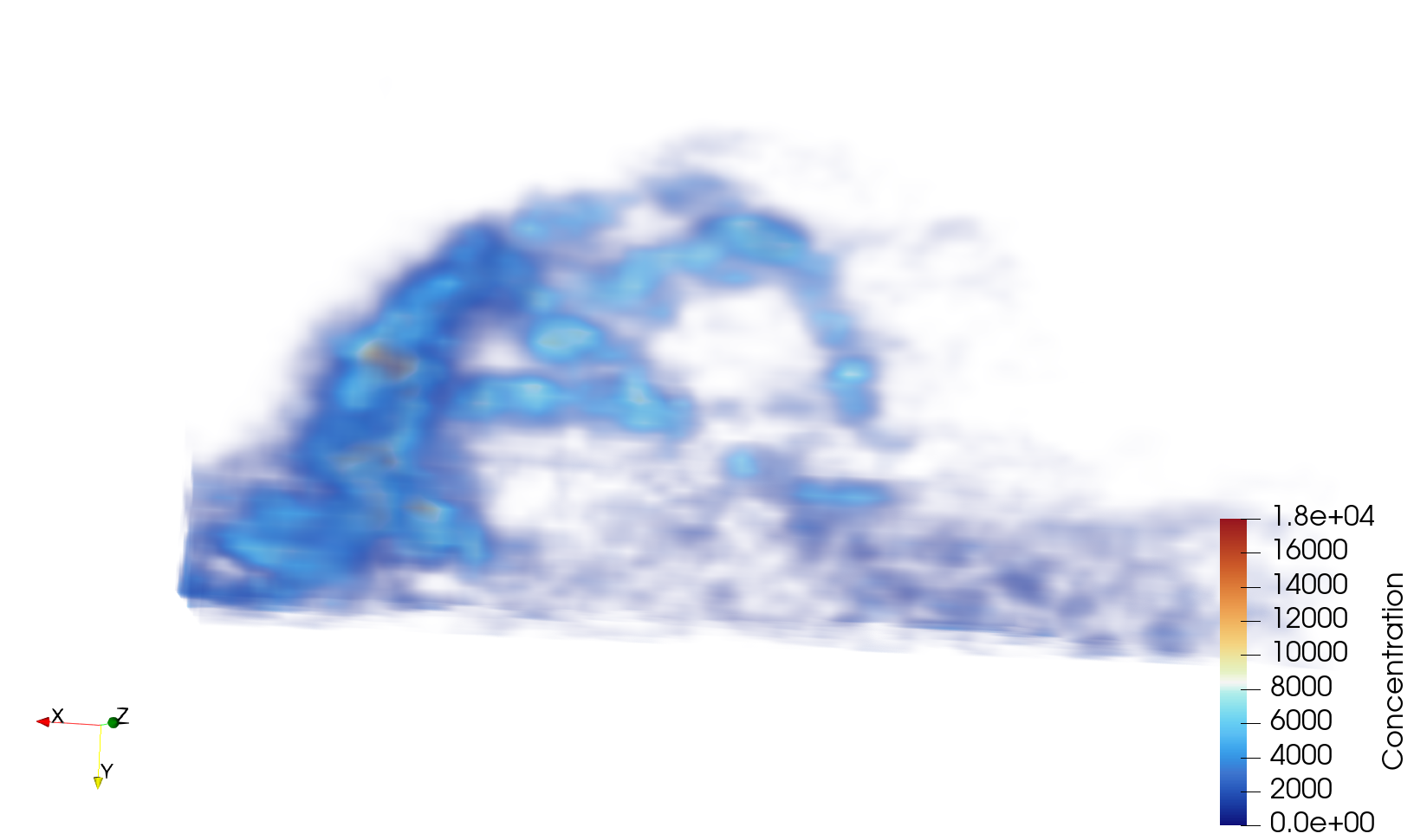}
    \end{subfigure}
    \hfill
    \begin{subfigure}[b]{0.48\textwidth}
        \centering
        \includegraphics[width=\textwidth]{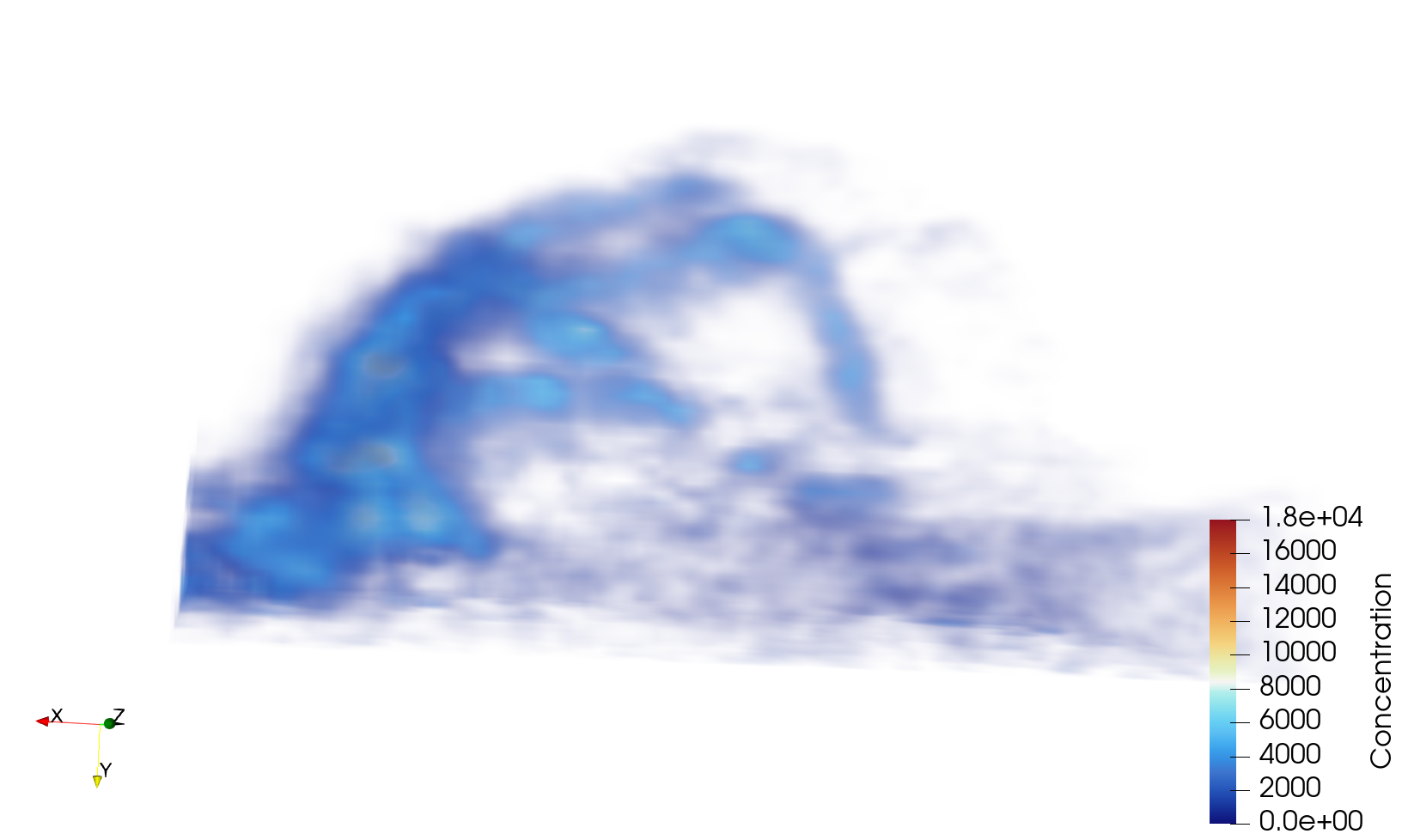}
    \end{subfigure}
    \caption{Visualisation of 3D ultrasound data at time $t = 12s$ before (left) and after (right) smoothing.}
    \label{fig:ultrasound}
\end{figure}

The entire dataset was smoothed using a Savitzky-Golay finite impulse response filter \cite{golayFilter}. 
This filter effectively reduces noise by applying a convolution process, in which consecutive subsets of adjacent 
data points are fitted with a low-degree polynomial via linear least squares. As a result, the data is smoothed and its 
precision enhanced, while preserving the underlying signal trends without distortion \cite{golayFilter}.
The difference before and after smoothing the concentration of tracer running through a mouse tumour at one time step 
is shown in Figure \ref{fig:ultrasound}.

After smoothing, the 3D data is projected onto a 2D domain by summing the discrete concentrations $C$ over the $z$-axis
\begin{equation}
    C_{2D}(x,y,t) = \sum_{i\in [0,77]} C(x,y,t_i).
\end{equation}
We chose to project onto the $x-y$-plane because it retains more information about the concentration flow compared to 
the other two planes. We can clearly observe the inflow, the transport through the organ, and the outflow of 
the tracer. 

Furthermore, we rescale the concentration from around $1.5\cdot10^{5}$ to
a maximum of $5$, to be  comparable to the synthetic data.
To ensure that there are no concentration values near the boundaries at the initial stage we add zero-padding with 16 cells, to reduce influence of boundary conditions by prohibiting inflow through the boundaries, i.e., it is assumed that the contrast agent bolus is fully inside the computational domain throughout simulation and reconstruction. To thoroughly treat inflow/outflow boundary conditions would require additional estimation of inflow/outflow location and velocities, and will be investigated in future work. With this padding, the domain size increases to $\Omega = 35.4mm \times 22.8mm$, 
and the number of cells to $118 \times 76$.

\begin{figure}
	\centering
	\includegraphics[width=0.5\textwidth]{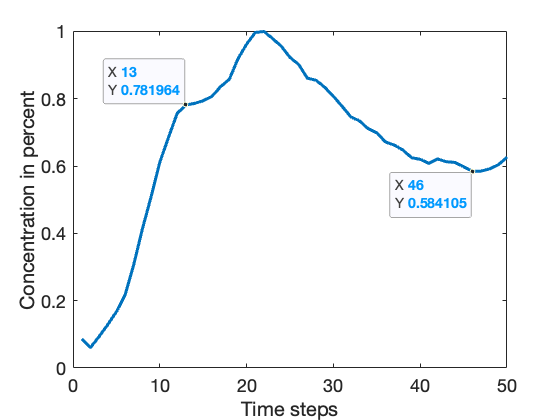}
	\caption{Percentage of concentration over time for the whole ultrasound measurement. Between measurement frame $13$ and $46$, indicated in the plot, 80\% of the total concentration are inside the computational domain.}
	\label{fig:concOverTime}
\end{figure}

Figure \ref{fig:concOverTime} illustrates the overall concentration over the time steps, clearly showing an inflow phase in the 
first half and an outflow phase in the second half. 
To ensure that the contrast agent is fully contained in the domain at initial time, we discard the initial and final frames, and start at time step $t_{s} = 13$ ($t = 6.96 \, s$),
and end at time step $t_{s} = 46$, ($t = 24.65 \, s$), such that 33 frames remain from the measured data. In this interval, at initial time approximately 80\% of the total concentration is inside the computational domain; experiments with setting this threshold to 70\% and 90\% yield similar results.

While using the spatial resolution of the measurement data for space discretization within the numerical solver is sufficient to be result in a stable simulation, the temporal resolution of the data is not high enough. For time discretization we use $120$ time steps, and linearly interpolate the measurement data. The remaining time interval is thus $(0, 17.64s)$, with a time step of $0.147 \, s$. 

The reconstructions were initialized with the parameters displayed in Table \ref{table:ultrasoundDataStarting}. For each model, we chose the optimizer that yielded the best results; for the two-compartment model this was steepest descent, while for advection-diffusion the Dai-Yuan conjugate gradient method was used.

\begin{table}[!ht]
    \centering
    \begin{tabular}{ l  c  c c  c  c  c  c } 
        \toprule
        model & $V^1 = (V^1_x \  V^1_y)$ & $V^2 = (V^2_x \  V^2_y)$ & $\kappa \vee D$ & $\{\lambda_i, i=1,2,3\}$& $\text{tol}_1$ & $\text{tol}_2$ & $\text{tol}_{\alpha}$\\ 
        \midrule
        Two-Compartment     & $\begin{pmatrix} 0.1 & -0.1 \end{pmatrix}$ & $\begin{pmatrix} 0.1 & 0.1 \end{pmatrix}$ & $12$   & $\{10^{-4},10^{-4},10^{-5}\}$ & $10^{-5}$ & $10^{-5}$ & $5\cdot10^{-7}$ \\  
        Advection-Diffusion & $\begin{pmatrix} 0.1 & 0.1 \end{pmatrix}$  & -                                         & $0.05$ & $ \{10^{-4},10^{-5}\}$      & $10^{-5}$ & $10^{-5}$ & $5\cdot 10^{-7}$ \\  
        \bottomrule
    \end{tabular}
    \caption{Initialization and parameters for the reconstructions of the ultrasound measurement data.}  
    \label{table:ultrasoundDataStarting}
\end{table}

\subsubsection{Results two-compartment model.}

\begin{figure}[!ht]
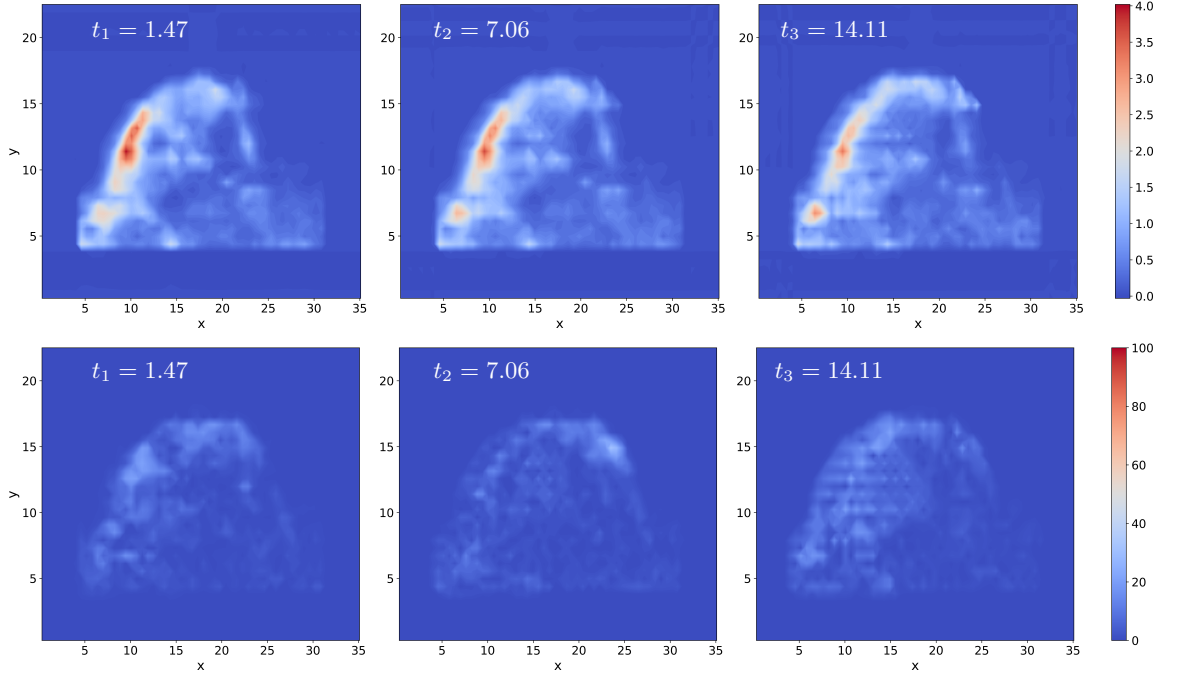

    \centering
    \begin{subfigure}[c]{\textwidth}
        \centering
        \begin{tikzpicture}

          \node[anchor=south west, inner sep=0] (img) 
          {\includegraphics[width=\textwidth]{Results/UltrasoundData/Ultrasound80/concentrations.png}};

          \begin{scope}[x={(img.south east)}, y={(img.north west)}]
          \node[font=\bfseries\small, text=white]
            at (0.12,0.9) {$t_1 = 1.47$};

          \node[font=\bfseries\small, text=white]
            at (0.415,0.9) {$t_2 = 7.06$};

          \node[font=\bfseries\small, text=white]
            at (0.72,0.9) {$t_3 = 14.11$};

          \end{scope}
        \end{tikzpicture}
    \end{subfigure}
        \begin{subfigure}[c]{\textwidth}
    	\centering
    	\begin{tikzpicture}
    		
    		\node[anchor=south west, inner sep=0] (img) 
    		{\includegraphics[width=\textwidth]{Results/UltrasoundData/Ultrasound80/concentrationRelErr.png}};

    		\begin{scope}[x={(img.south east)}, y={(img.north west)}]
    			
    			\node[font=\bfseries\small, text=white]
    			at (0.12,0.9) {$t_1 = 1.47$};
    			
    			\node[font=\bfseries\small, text=white]
    			at (0.415,0.9) {$t_2 = 7.06$};
    			
    			\node[font=\bfseries\small, text=white]
    			at (0.715,0.9) {$t_3 = 14.11$};
    			
    		\end{scope}
    	\end{tikzpicture}
    \end{subfigure}
    \caption{Two-compartment model: reconstructed concentrations (top row) with space-dependent $\kappa$  and relative error (in percent) between the final simulated and measured concentrations (bottom row).
    }
    \label{fig:ultraConcentration}
\end{figure}

As the conversion domain is typically unknown, we only reconstruct space-dependent velocities and transfer coefficient.  The algorithm terminated after four outer iterations, as no feasible step size was found. During optimization, cost functional values and gradient norms decreased significantly, see Table~\ref{table:resultUS}, indicating a successful reconstruction, even though the final result is still far from a stationary point.  To evaluate the reconstruction, simulated concentrations using the identified parameter fields are shown in Figure \ref{fig:ultraConcentration}, together with the relative errors to the measured concentrations.
Even though the mismatch is quite high, the main flow of the tracer agent has 
been reconstructed. The high errors can be explained by the mismatch in overall concentration, given that the reconstructed 
problem has around $20\%$ less tracer concentration because it cannot model the inflow over time that occurs in the data.

\begin{figure}[!htp]
	\centering
	\includegraphics[width=\textwidth]{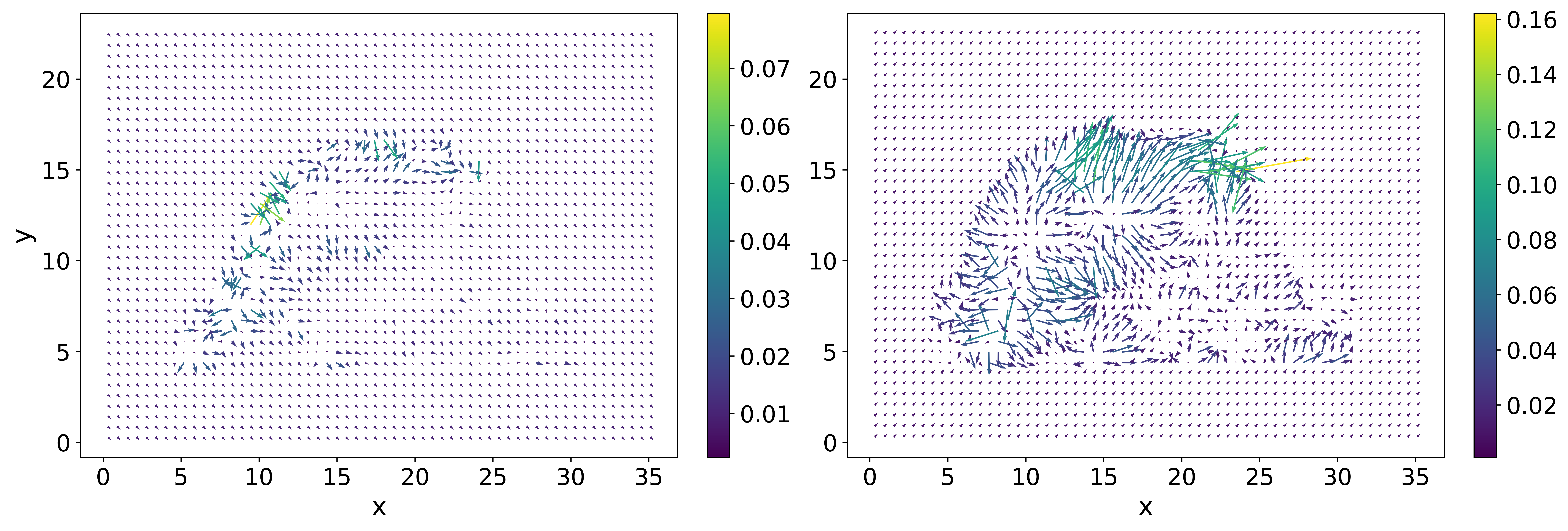}
	\caption{Two-compartment model: reconstructed arterial (left) and venous (right) velocities. Arrows (one per cell) indicate the direction, coloring the magnitude. Nonzero velocities are reconstructed only in areas with contrast agent concentrations.}
	\label{fig:QuiverUS_k}
\end{figure}

\begin{figure}[!ht]
    \centering
    \includegraphics[width=0.65\textwidth]{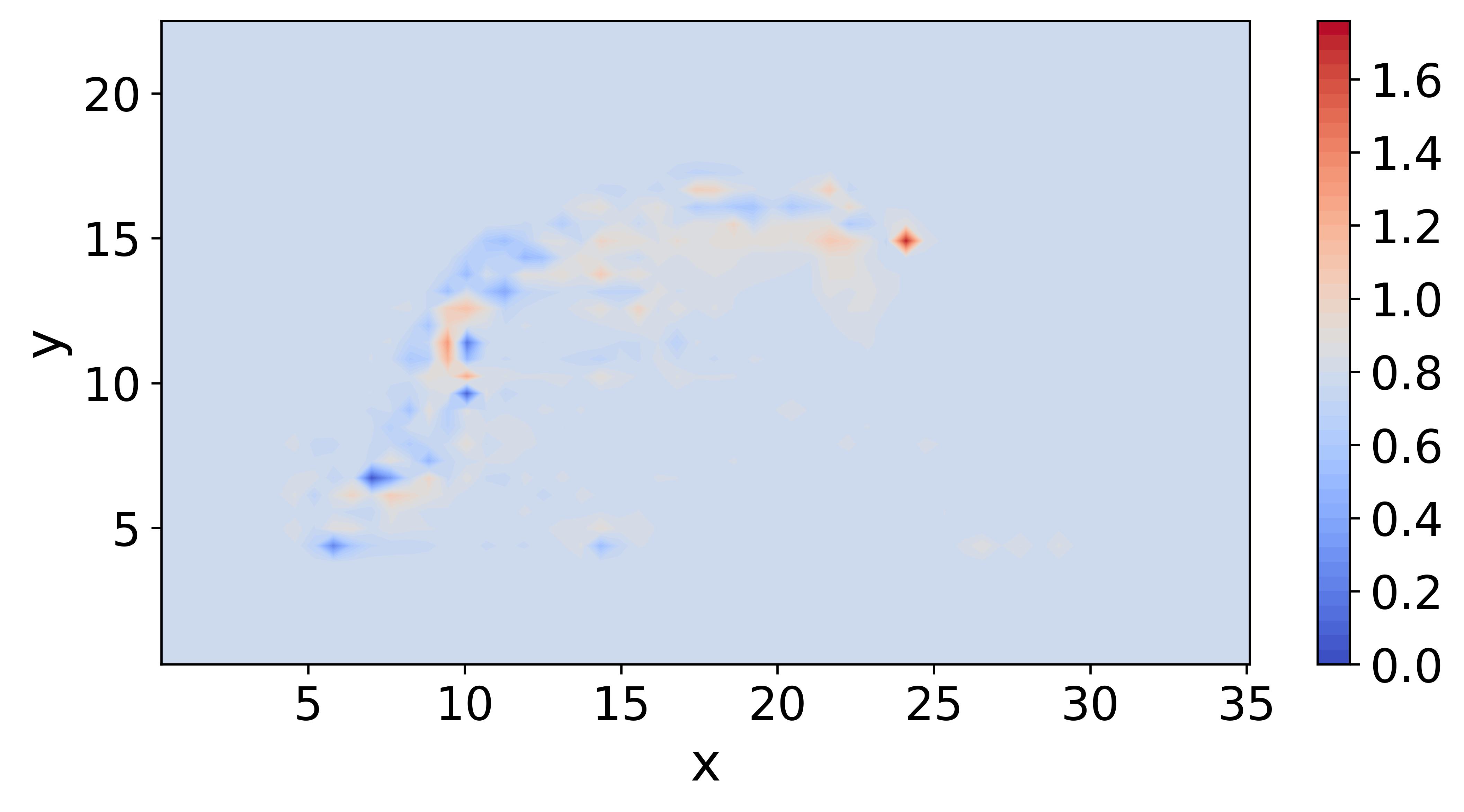}
    \caption{Two-compartment model: reconstructed space-dependent transfer coefficient $\kappa$.}
    \label{fig:kappaRec}
\end{figure}

The reconstructed velocities and transfer coefficient are 
shown in Figures~\ref{fig:QuiverUS_k} and~\ref{fig:kappaRec}. As expected, we observe higher velocity activity in the arterial component at 
the beginning of the domain and lower activity towards the end. The opposite is true for the velocity in the venous component. 
However, velocity magnitudes are not entirely as physiologically expected: venous velocities should be lower than arterial ones. The deviation in the reconstruction might be due to the projection of 3D data to a 2D plane, as well as the relatively short measurement time span used for the reconstruction.
The transfer coefficient shows higher values towards the center of the organ and lower values
at the beginning, which matches the expected behavior.

In conclusion, the reconstructions derived from ultrasound data appear promising, delivering qualitatively reasonable 
 velocity fields and transfer coefficient. 

\subsubsection{Results advection-diffusion model.}

\begin{figure}[!ht]
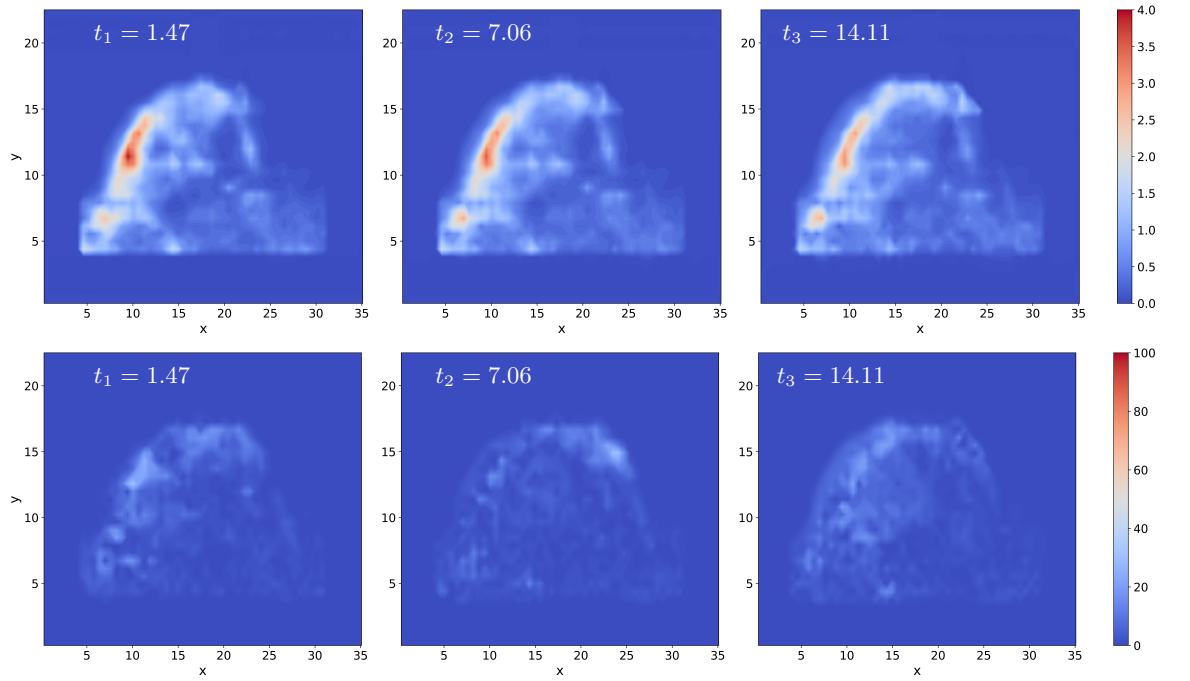

	\centering
	\begin{subfigure}[c]{\textwidth}
		\centering
		\begin{tikzpicture}
			
			\node[anchor=south west, inner sep=0] (img) 
			{\includegraphics[width=\textwidth]{Results/UltrasoundData/Ultrasound80/concentrationsAD.png}};
	
			\begin{scope}[x={(img.south east)}, y={(img.north west)}]
				\node[font=\bfseries\small, text=white]
				at (0.12,0.9) {$t_1 = 1.47$};
				
				\node[font=\bfseries\small, text=white]
				at (0.415,0.9) {$t_2 = 7.06$};
				
				\node[font=\bfseries\small, text=white]
				at (0.72,0.9) {$t_3 = 14.11$};
				
			\end{scope}
		\end{tikzpicture}
	\end{subfigure}
	\begin{subfigure}[c]{\textwidth}
		\centering
		\begin{tikzpicture}
			
			\node[anchor=south west, inner sep=0] (img) 
			{\includegraphics[width=\textwidth]{Results/UltrasoundData/Ultrasound80/concentrationRelErrAD.png}};

			\begin{scope}[x={(img.south east)}, y={(img.north west)}]
				
				\node[font=\bfseries\small, text=white]
				at (0.12,0.9) {$t_1 = 1.47$};
				
				\node[font=\bfseries\small, text=white]
				at (0.415,0.9) {$t_2 = 7.06$};
				
				\node[font=\bfseries\small, text=white]
				at (0.715,0.9) {$t_3 = 14.11$};
				
			\end{scope}
		\end{tikzpicture}
		\end{subfigure}
	\caption{Advection-diffusion model: reconstructed concentrations (top row) with space-dependent diffusivity and relative error (in percent) between the final simulated and measured concentrations (bottom row).}
	\label{fig:ultraConcentrationAD}
\end{figure}

\begin{figure}[!ht]
\centering
\begin{subfigure}[c]{0.49\textwidth}
\includegraphics[width=\textwidth]{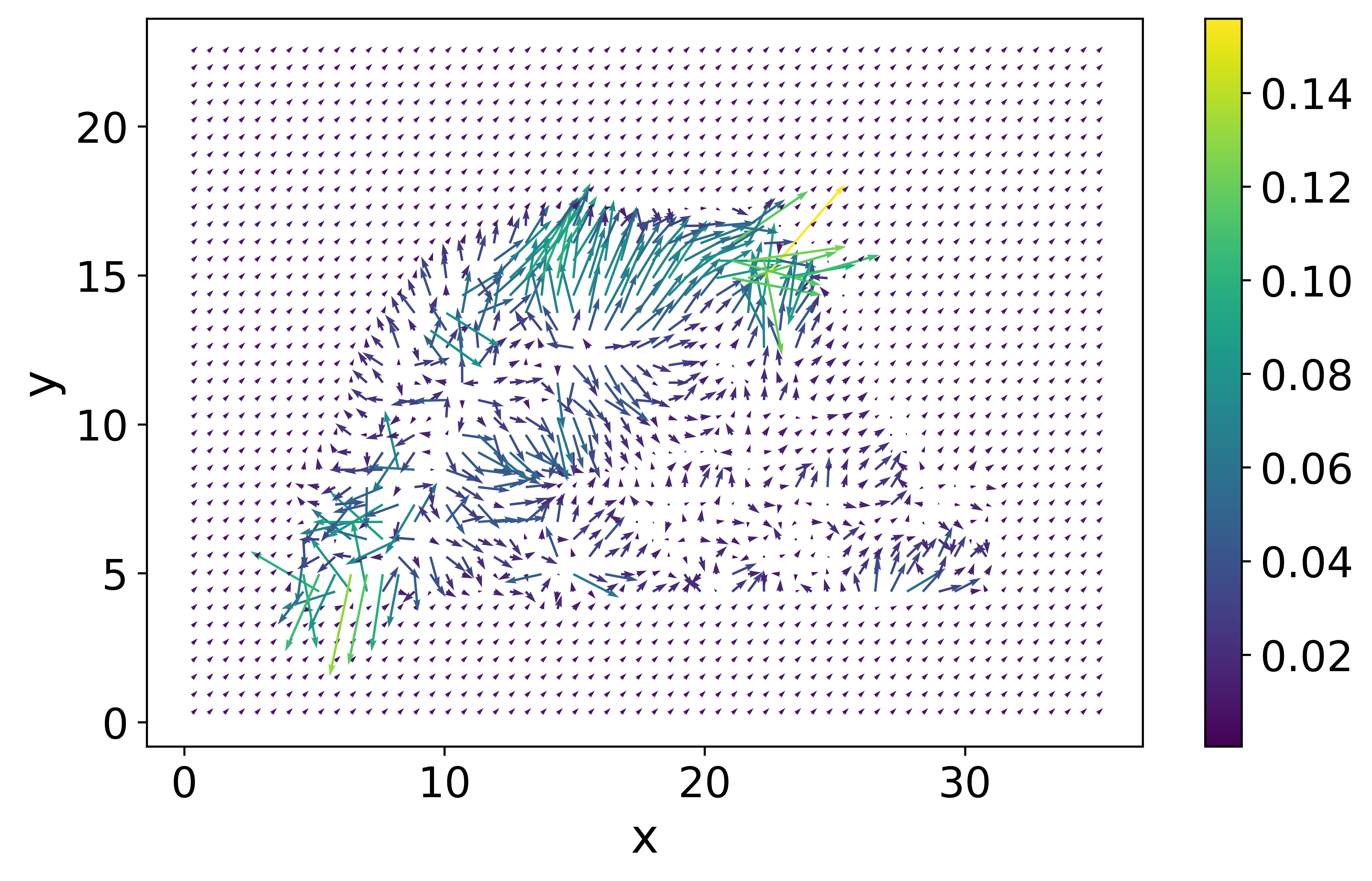}
\end{subfigure}
\begin{subfigure}[c]{0.49\textwidth}
	\includegraphics[width=\textwidth]{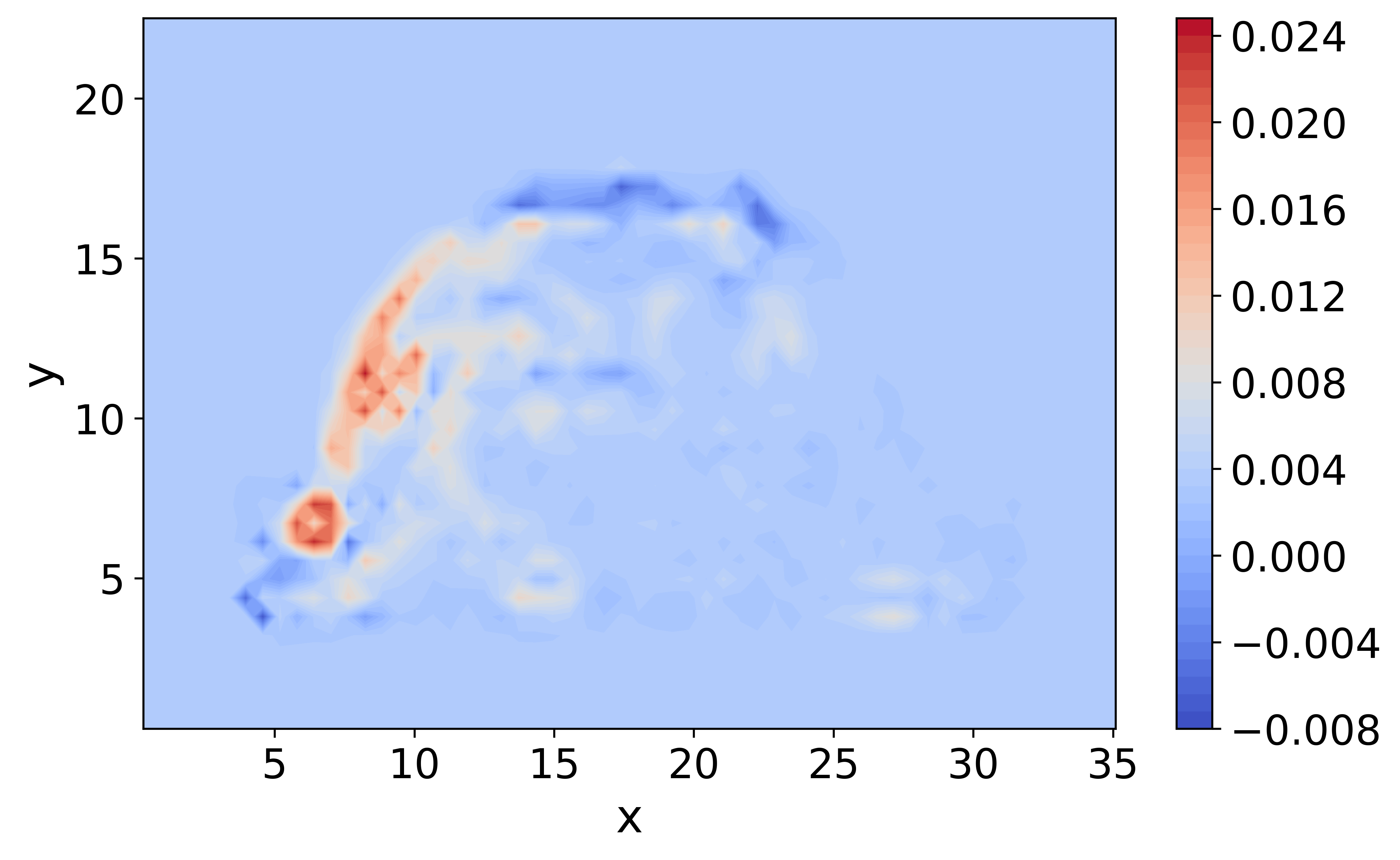}
\end{subfigure}
\caption{Advection-diffusion model: reconstructed velocity (left) and diffusivity (right). }
	\label{fig:ultraVelAndDiffAD}
\end{figure}

For comparison we use the advection-diffusion model for reconstruction of spatially varying velocity and diffusivity.
The intial parameters can be found in Table \ref{table:ultrasoundDataStarting}. After 20 outer iterations, cost functional and the gradient norm of the velocity show reasonable decrease, 
compare Table \ref{table:resultUS}. The gradient norm in diffusivity is still quite high but was decreased by roughly $75\%$.
The resulting concentrations are shown in Figure \ref{fig:ultraConcentrationAD} together with the relative error to the measured concentrations. The movement of contrast agent and 
the mismatch are similar to the concentrations reconstructed with the two-compartment model. Reconstructed velocity and diffusivity are shown in  Figure~ \ref{fig:ultraVelAndDiffAD}. While the velocity is again in a physiologically meaningful range, the identified diffusivity is negative in some parts. As this is physically not correct, it indicates that the two-compartment model is better suited for the identification of perfusion parameters. 


\begin{table}[!ht]
    \centering
    \begin{tabular}{  p{2.8cm}  p{2cm}  p{2cm}  p{3cm}  p{0.8cm} p{0.8cm} } 
        \toprule
         & $||\nabla j_{v}||_{\text{start}}$ & $||\nabla j_{v}||_{\text{final}}$ & $||\nabla j_{\kappa}||$ or $||\nabla j_{D}||$ & $J_{\text{start}}$ & $J_{\text{final}}$ \\ 
        \midrule
        Two-compartment  & $19.94$ & $0.204$ & $0.002$ & $14.53$ & $8.29$ \\
        Advection-diffusion & $9.67$  & $0.19$ & $0.128$ & $12.42$ & $6.96$  \\
        \bottomrule
    \end{tabular}
    \caption{Initial and final gradient norms and objective function values for the two-compartment and advection-diffusion models applied to measurement data, indicating a reasonable reconstruction in both cases.}  
    \label{table:resultUS}
\end{table}

\section{Conclusion}

We have derived adjoint equations for efficient gradient computation for two models for tracer transport in blood flow: a simple, more commonly used advection-diffusion model, as well as a more realistic two-compartment model, that, to the best of our knowledge, was not used for parameter identification in ultrasound before. We demonstrated feasibility for reconstructing flow parameters from contrast agent concentrations for synthetic as well as measurement data. For the latter, both models, the advection-diffusion and two-compartment
model, are able to reconstruct perfusion-related parameters. 
While the qualitative behavior of the reconstructed arterial and venous velocities of the two-compartment model matches the expected behavior, an extension of the implementation to 3D as well as additional identification of inflow and outflow locations and velocities is needed to improve reconstructed velocity magnitudes.
However, the two-compartment model yields a reasonable transfer coefficient function, in contrast to the partially negative diffusivity reconstructed by the advection-diffusion model. This indicates that the two-compartment model is more suitable for use in medical applications. However, to  confirm the relation between the identified transfer coefficient and actual tumor perfusion,  the model needs to be validated against, e.g., histology data. Together with the extension of the implementation to 3D in space and a Bayesian approach for uncertainty quantification, this is ongoing work.


%
%


\paragraph{Funding.}{
Research reported in this publication
was supported by the National Cancer Institute of the National Institutes of Health under Award Number
R01CA286505. The content is solely the responsibility of the authors and does not necessarily represent the
official views of the National Institutes of Health.}




\printbibliography

\appendix

\section{Derivation of Adjoint System and Gradient for the Advection-Diffusion Model}\label{sec:adjointAdvecDiff}
We follow the formal Lagrange principle~\cite{troltzsch}, and define the Lagrangian $L(u,c,p)$ 
\begin{equation}
  \begin{split}
    L(u,c,p) &= J(u,c) - (p, u_{t} + \nabla \cdot (V \, u) - \nabla \cdot (D \, \nabla u))_{L_{2}(Q)}\\
    J(u, (V,D)) &= \frac{1}{2} ||u - C_\text{meas}^\delta||_{L_{2}(Q)}^{2} + \frac{\lambda_{1}}{2} ||V||_{L_{2}(\Omega)^d}^{2} + \frac{\lambda_{2}}{2}||D||_{L_{2}(\Omega)}^{2},
  \end{split}
\end{equation}
where $p$ is the adjoint state. Here, the problem is posed on the space-time cylinder $Q=\Omega \times (0,T)$. Its boundary is denoted by $\Sigma = \partial \Omega \times (0,T)$ . We assume that $V$,  $D$ and the resulting solution $u$ are sufficiently regular,  
and use the standard $L_2$ scalar product given by $(u,v)_{L_{2}(Q)} = \int_{Q} u(x)  v(x) \,dx$. To derive the adjoint equation  we
compute the derivative of $L(u,c,p)$ with respect to the state $u$ in direction $d$:

\begin{equation}
  \begin{split}
    L_{u}(u,c,p) \, d &= J_{u}(u,c) \, d - (p, d_{t} + \nabla \cdot (V \, d) - \nabla \cdot (D \, \nabla d))_{L_{2}(Q)}\\
                 &= J_{u}(u,c) \, d - \int_{0}^{T} \int_{\Omega} p \, (d_{t} + \nabla \cdot (V \, d) -  \nabla \cdot (D \, \nabla d))\,dx\,dt,\\
  \end{split}
\end{equation}
Using  integration by parts  in time we get 
\begin{equation}
  \begin{split}
    \int_{0}^{T}\int_{\Omega} p \, d_{t}\,dx\,dt = &\int_{\Omega}[d(x,t)\, p(x,t)]_{0}^{T} \,dx - \int_{0}^{T}\int_{\Omega} d \, p_{t}\,dx\,dt\\
                                                                = &\int_{\Omega}(d(x,T)\, p(x,T) -d(x,0)\, p(x,0))\,dx \\
                                                                  &- \int_{0}^{T}\int_{\Omega} d \, p_{t}\,dx\,dt.\\
  \end{split}
\end{equation}
Applying Green's identity \cite{troltzsch} and  integration by parts, we obtain, with $\nu$ denoting the outer unit normal vector on $\Sigma$,  
\begin{equation}
  \begin{split}
    \int_{0}^{T}\int_{\Omega} p\, \nabla \cdot (V \, d) \,dx\,dt = &- \int_{0}^{T} \int_{\Omega} \nabla p \cdot V \, d \,dx\,dt\\
                                                                                 &+ \int_{0}^{T} \int_{\Sigma} p \, V \cdot \nu \, d \,dS\,dt\\
    \int_{0}^{T}\int_{\Omega} p \, \nabla \cdot (D \, \nabla d) \,dx\,dt = &- \int_{0}^{T} \int_{\Omega} \nabla p \cdot (\nabla d \, D) \,dx\,dt\\
                                                                                         &+ \int_{0}^{T} \int_{\Sigma} p \, D \, \partial_{\nu} d \,dS\,dt\\
                                                                                       = &+ \int_{0}^{T} \int_{\Omega} \nabla \cdot (D \, \nabla p) \, d \,dx\,dt\\
                                                                                         &- \int_{0}^{T} \int_{\Sigma} \partial_{\nu} p \, D \, d \,dS\,dt\\
                                                                                         &+ \int_{0}^{T} \int_{\Sigma} p \, D \, \partial_{\nu} d \,dS\,dt\\
  \end{split}
\end{equation}
Combined with the derivative of the cost functional with respect to $u$ we get
\begin{equation}
  \begin{split}
    L_{u}(u,c,p)\, d = \,&(u - C_\text{meas}^{\delta}, d)_{L_{2}(Q)} + (p_t, d)_{L_{2}(Q)} + (\nabla p \cdot V, d)_{L_{2}(Q)}\\ 
                          &+ (\nabla \cdot (D \, \nabla p), d)_{L_{2}(Q)} \\
                          &- \int_{\Sigma} d \, (p \, V \cdot \nu + \partial_{\nu} p \, D) - p \, D \, \partial_{\nu} d \,dS\,dt\\
                          &- \int_{\Omega}(d(x,T)\, p(x,T) -d(x,0)\, p(x,0))\,dx.
  \end{split}
\end{equation}
To identify the strong formulation of the adjoint equation, we 
first choose $d \in C^{\infty}_{0}(\Omega)$, which yields $d = \partial_{\nu} d = 0$ on $\Gamma$ and $d(T) = d(0) = 0$ on $\Omega$.
        \begin{equation}
          \begin{split}
            &\int_{0}^{T}\int_{\Omega} d \, [(u-u^{\delta}) + (p_t + \nabla p \cdot V + \nabla \cdot (D \, \nabla p))] \,dx\,dt = 0\\
            &\Rightarrow -p_t - \nabla p \cdot V - \nabla \cdot (D \, \nabla p) = u - u^{\delta} \text{ on } Q,
          \end{split}
        \end{equation}
as $C^{\infty}_{0}(\Omega)$ lies dense in $L_2(\Omega)$.
 Next, we allow $d$ to vary on the boundary, and get 
        \begin{equation}
          p \, V \cdot \nu + \partial_{\nu}p \, D = 0 \text{ on } \Gamma.
        \end{equation}
Including the terminal condition $p(\cdot, T) = 0$ we arrive at the adjoint equation
\begin{equation}
  \begin{split}
    -p_t - \nabla p \cdot V - \nabla \cdot (D \, \nabla p) &= u - u^{\delta} \text{ on } Q\\
    p(\cdot,T) &= 0 \\
    p \, V \cdot \nu + \partial_{\nu}p \, D &= 0 \text{ on } \Gamma \\
  \end{split}
\end{equation}
For the implementation, the backward-in-time adjoint equation is re-writte using the time transform 
$\tau = T - t$, where $T$ is the final time and $t \in [0,T]$, to allow using the same solver for state and adjoint equations.
This yields
\begin{equation}
    p_{\tau} - \nabla p \cdot V - \nabla \cdot (D \, \nabla p) = (u(\tau) - C_\text{meas}^{\delta}(\tau)),
\end{equation}
with initial condition $p(\cdot, 0) = 0$.
Next we derive the descent direction for the optimization. For this we calculate the derivative of the Lagrangian with respect
to the parameters $c = (V, D)$ in direction $d = (d_V, d_D)$. 
\begin{equation}
	\begin{split}
		L_{V}(u,p,c)\, d_{V} &= J_{V}(u,c)\, d_{V} - \left(p, \nabla \cdot \left(d_V\, u \right)\right)_{L_2(Q)}  \\
		&= J_{V}(u,c)\, d_{V} + \int_0^T\int_{\Omega} \frac{\partial p}{\partial x} \, d_{V} \, u \,dx\,dt - \int_{\Sigma} p \, d_{V} \, \nu \, u \,dS\,dt\\
		&= J_{D}(u,c)\, d_{D} - \int_0^T \int_{\Omega} \nabla p \cdot \nabla u \, d_D \,dx\,dt + \int_{\Sigma} p \, \partial_{\nu} u \, d_D \,dS\,dt.
	\end{split}
\end{equation}
%
Using further that the derivatives of the cost functional are given by
\begin{equation}
  \begin{split}
    J_{V}(u,c)\, d_{V} &= (\lambda \, V, d_{V})_{L_2(\Omega)}\\
    J_{D}(u,c)\, d_{D} &= (\lambda_2 \, D, d_{D})_{L_2(\Omega)}.
  \end{split}
\end{equation}
We identify the gradient of the reduced cost functional by Riesz representation theorem as
\begin{equation}
  \begin{split}
    \nabla j(V) &= \int_0^T \lambda_1^{T} \, V + \nabla p \, u  \,dt\\
    \nabla j(D) &= \int_0^T \lambda_2 \, D - \nabla p \cdot \nabla u \,dt.
  \end{split}
\end{equation}

\section{Derivation of adjoint system and gradient for the two-compartment model}\label{sec:adjoint2C}
Using the cost functional we define the Lagrangian
\begin{equation}
  \begin{split}
    L((u,w),(p,q),c) &= J(u,w,c) - (p, u_{t} + \nabla \cdot (v^1 \, u) + \kappa \, u)_{L_{2}(Q)}\\
                     &\hspace*{1.9cm}           - (q, w_{t} + \nabla \cdot (v^2 \, w) - \kappa \, u)_{L_{2}(Q)}\\
    J(u,w,(v^1,v^2,\kappa)) &= \frac{1}{2} ||u + w - C_\text{meas}^{\delta}||_{L_{2}(Q)}^{2} 
                                + \frac{\lambda_{1}}{2} ||V^1||_{L_{2}(\Omega)^d}^{2} + \frac{\lambda_{2}}{2}||V^2||_{L_{2}(\Omega)^d}^{2}
                                + \frac{\lambda_{3}}{2} ||\kappa||_{L_{2}(\Omega)}^{2},
  \end{split}
\end{equation}
with $p$ and $q$ the adjoint states and $C_\text{meas}^{\delta}$ are the measurements over time. 
To derive the adjoint equation, we again assume sufficient regularity of all involved functions and consider the derivatives of $L(u,w,c,p)$ with respect to the states $u$ and $w$, in a direction $d$ respectively. 
\begin{equation}
  \begin{split}
    L_{u}((u,w),(p,q))\, d = J_{u}(u,w,c) \, d &- (p, d_{t} + \nabla \cdot (v^1 \, d) + \kappa \, d)_{L_{2}(Q)} \\
                                &- (q, - \kappa \, d)_{L_{2}(Q)}\\
  \end{split}
\end{equation}
Using  integration by parts we get
\begin{equation}
  \begin{split}
     \int_{0}^{T}\int_{\Omega} p \, d_{t}\,dx\,dt =  &\int_{\Omega}[d(x,t)\, p(x,t)]_{0}^{T} dx - \int_{0}^{T}\int_{\Omega} d \, p_{t}\,dx\,dt\\
                                                                =  &\int_{\Omega}(d(x,T)\, p(x,T) -d(x,0)\, p(x,0))\,dx\\ 
                                                                   &- \int_{0}^{T}\int_{\Omega} d \, p_{t}\,dx\,dt.\\
  \end{split}
\end{equation}
Applying Green's identity we further obtain
\begin{equation}
     \int_{0}^{T}\int_{\Omega} p \, \nabla \cdot (V^1 \, d) \,dx\,dt = - \int_{0}^{T} \int_{\Omega} \nabla p \cdot V^1 \, d \,dx\,dt + \int_{0}^{T} \int_{\Sigma} p \, V^1 \cdot \nu \, d \,dS\,dt
\end{equation}
With the derivative of the cost functional with respect to $u$ in direction $d$,
\begin{equation}
  \begin{split}
    J_u(u,w,c) \, d &= (u + w - C_\text{meas}^{\delta},d)_{L_{2}(Q)},
  \end{split}
\end{equation}
this yields
\begin{equation}
  \begin{split}
    L_{u}((u,w),(p,q),c)\, d = & \hspace*{0.1cm}(u + w - C_\text{meas}^{\delta}, d)_{L_{2}(Q)} + (p_t, d)_{L_{2}(Q)}\\ 
                                  &+ (\nabla p \cdot V^1, d)_{L_{2}(Q)} + (\kappa \, (q - p), d)_{L_{2}(Q)} \\
                                  &- \int_{\Sigma} d\, p \, V^1 \cdot \nu \,dS\,dt\\
                                  &- \int_{\Omega} (d(x,T)\, p(x,T) -d(x,0)\, p(x,0)) \,dx.    
  \end{split}
\end{equation}
Doing the same for the second state component $w$ yields
\begin{equation}
  \begin{split}
    L_{w}((u,w),(p,q),c)\, d = & \hspace*{0.1cm}(u + w - C_\text{meas}^{\delta}, d)_{L_{2}(Q)} + (q_t, d)_{L_{2}(Q)} + (\nabla q \cdot V^2, d)_{L_{2}(Q)}\\
                                  &- \int_{\Sigma} d\, q \, V^2 \cdot \nu \,dS\,dt\\
                                  &- \int_{\Omega} (d(x,T)\, q(x,T) -d(x,0)\, q(x,0)) \,dx. 
  \end{split}
\end{equation}
Choosing  $d \in C^{\infty}_{0}(Q)$  yields $d = \partial_{\nu} d = 0$ on $\Gamma$ and $d(T) = d(0) = 0$ on $\Omega$.
    \begin{equation}
      \begin{split}
        &\int_{0}^{T}\int_{\Omega} d \, [(u + w -C_\text{meas}^{\delta}) + (p_t + \nabla p \cdot V^1 + \kappa \, (q - p))] \,dx\,dt = 0\\
        &\Rightarrow -p_t - \nabla p \cdot V^1 - \kappa \, (q - p) = u + w - C_\text{meas}^{\delta} \text{ on } Q\\
        &\int_{0}^{T}\int_{\Omega} d \, [(u+w-C_\text{meas}^{\delta}) + (q_t + \nabla q \cdot V^2)] \,dx\,dt = 0\\
        &\Rightarrow -q_t - \nabla q \cdot V^2 = u+w - C_\text{meas}^{\delta} \text{ on } 
      \end{split}
    \end{equation}
  For the boundary conditions we get, due to $\partial_{\nu} u|_{\Gamma} = 0$, we know that 
        $\partial_{\nu} d|_{\Gamma} = 0$, and for $d|_{\Gamma}$  we get 
        \begin{equation}
          \begin{split}
            p \, V^1 \cdot \nu &= 0 \Rightarrow p |_{\Gamma} = 0\\
            q \, V^2 \cdot \nu &= 0 \Rightarrow q |_{\Gamma} = 0.
          \end{split}
        \end{equation}
Combined this yields the strong form of the adjoint equation
\begin{equation}
  \begin{split}
    -p_t - \nabla p \cdot V^1 + \kappa \, p &= u + w - C_\text{meas}^{\delta} + \kappa q\\
    -q_t - \nabla q \cdot V^2 &= u+w - C_\text{meas}^{\delta}\\
    p(x,T) &= 0 \\
    q(x,T) &= 0, \\
  \end{split}
\end{equation}
together with homogeneous Dirichlet boundary conditions.
Similar to the advection-diffusion model, we use the time transform $\tau = T - t$, where $T$ is the final time and $t \in [0,T] $ to obtain
\begin{equation}
  \begin{split}
    p_{\tau} - \nabla p \cdot V^1 + \kappa \, p &= (u(\tau) + w(\tau)- C_\text{meas}^{\delta}(\tau)) + \kappa q\\
    q_{\tau} - \nabla q \cdot V^2 &= (u(\tau) + w(\tau) - C_\text{meas}^{\delta}(\tau)).
  \end{split}
\end{equation}
To ensure that we can use the same solver as for the state equation, we further to reorder the equation as
\begin{equation}
  \begin{split}
    p_{\tau} - \nabla \cdot (v^1 \, p) + (\kappa + \nabla \cdot v^1)\, p &= (u(\tau) + w(\tau) - C_\text{meas}^{\delta}(\tau)) + \kappa \, q\\
    q_{\tau} - \nabla \cdot (v^2 \, q) + \nabla \cdot v^2 \cdot q &= (u(\tau) + w(\tau) - C_\text{meas}^{\delta}(\tau)).
  \end{split}
\end{equation}
After deriving the adjoint problem we can deduce the descent direction needed for the
optimization by calculating the derivative of the Lagrangian with respect
to $c$. Because $c$ contains three variables ($V^1,V^2, \kappa$) we calculate the partial derivatives
individually. The resulting expressions are
\begin{equation}
  \begin{split} 
    L_{V^1}(u, w,p,c)\, d_{V^1} &= J_{V^1}(u,w,c)\, d_{V^1} - \left(p, \nabla \cdot \left(d_{V^1} \, u \right)\right)_{L_2(Q)} \\
                                    &= J_{V^1}(u,w,c)\, d_{V^1} + \int_0^T\int_{\Omega} \frac{\partial p}{\partial x} \, d_{V^1} \, u \,dx\,dt - \int_{\Sigma} p \, d_{V^1} \, \nu \, u \,dS\,dt\\
    L_{V^2}(u, w,q,c)\, d_{V^2} &= J_{V^2}(u,w,c)\, d_{V^2} - \left(q, \nabla \cdot \left(d_{V^2}\, w \right)\right)_{L_2(Q)}\\
                                    &= J_{V^2}(u,w,c)\, d_{V^2} + \int_0^T\int_{\Omega} \frac{\partial q}{\partial x} \, d_{V^2} \, w \,dx\,dt - \int_{\Sigma} q \, d_{V^2} \, \nu \, w \,dS\,dt\\
    L_{\kappa}(u, w,q,c) \, d_{\kappa} &= (\lambda_{3}\, \kappa, d_{\kappa})_{L_{2}(\Omega)} - (p,d_{\kappa} \, u)_{L_{2}(Q)} - (q,-d_{\kappa} \, u)_{L_{2}(Q)}.
  \end{split}
\end{equation}
From this we can identify the gradient of the reduced cost functional by the Riesz representation theorem as
\begin{equation}
  \begin{split}
    \nabla j(V^1) &= \int_0^T \lambda_1 \, V^1 + \nabla p \, u \,dt\\
    \nabla j(V^2) &= \int_0^T \lambda_2 \, V^2 + \nabla q \, w \,dt\\
    \nabla j(\kappa) &= \int_0^T \lambda_{3}\, \kappa + (q - p) \, u \,dt.
  \end{split}
\end{equation}

\end{document}